\documentclass[11pt,a4paper]{article}

\usepackage{mathptmx} 
\DeclareMathAlphabet{\mathcal}{OMS}{cmsy}{m}{n}
\usepackage[top=25mm,bottom=35mm,left=20mm,right=20mm]{geometry} 

\usepackage{graphicx}
\usepackage{amssymb,amsmath,amsfonts,amsthm}
\usepackage{bm}
\usepackage[english]{babel}
\usepackage{url}
\usepackage{doi}
\usepackage{tikz}
\usepackage{mathtools}
\mathtoolsset{showonlyrefs=true}

\newtheorem{prop}{Proposition}
\newtheorem{remark}{Remark}

\newcommand{\R}{\mathbb R}

\newcommand{\hl}{\ensuremath{\frac{1}{2}}}
\renewcommand{\d}{\ensuremath{\delta}}

\newcommand{\px}{\partial _x}
\newcommand{\bra}[1]{\left( #1\right)}
\newcommand{\dx}{\mathrm{d}x}
\newcommand{\hfline}{[-L,L]}
\newcommand{\fracdel}[2]{\frac{\delta
#1}{\delta #2}}
\newcommand{\dsum}{{\sum_{n=0}^N} \, ''}

\begin{document}

\title{Geometric numerical integrators for Hunter--Saxton-like equations}
\author{
Yuto Miyatake\thanks{Department of Computational Science and Engineering, 
              Graduate School of Engineering, 
              Nagoya University,  Furo-cho, Chikusa-ku, 
              464-8603 Nagoya, Japan, 
\href{mailto:miyatake@na.nuap.nagoya-u.ac.jp}{miyatake@na.nuap.nagoya-u.ac.jp}} ,\
David Cohen\thanks{Matematik och matematisk statistik, 
              Ume{\aa} universitet, 90187 Ume{\aa}, 
              Sweden, \href{mailto:david.cohen@umu.se}{david.cohen@umu.se}}
              \thanks{Department of Mathematics, University of Innsbruck, 
              6020 Innsbruck, Austria, \href{mailto:david.cohen@uibk.ac.at}{david.cohen@uibk.ac.at}}, \
Daisuke Furihata\thanks{Cybermedia Center,
              Osaka University, Machikaneyama 1-32, Toyonaka,
              Osaka, 560-0043, Japan,\href{mailto:furihata@cmc.osaka-u.ac.jp}{furihata@cmc.osaka-u.ac.jp}}
\ and 
Takayasu Matsuo\thanks{Department of Mathematical Informatics, 
              Graduate School of System of Information Science and Technology, 
              The University of Tokyo, 7-3-1 Hongo Bunkyo-ku, 
              113-8656 Tokyo, Japan, 
\href{mailto:matsuo@mist.i.u-tokyo.ac.jp}{matsuo@mist.i.u-tokyo.ac.jp}} }

\maketitle

\begin{abstract}
We present novel geometric numerical integrators for Hunter--Saxton-like equations 
by means of new multi-symplectic formulations and known Hamiltonian structures 
of the problems. 
We consider the Hunter--Saxton equation, the modified Hunter--Saxton equation, 
and the two-component Hunter--Saxton equation. 
Multi-symplectic discretisations based on these new formulations of the problems 
are exemplified by means of the explicit Euler box scheme, and Hamiltonian-preserving discretisations 
are exemplified by means of the discrete variational derivative method. 
We explain and justify the correct treatment of boundary conditions in a unified manner. 
This is necessary for a proper numerical implementation of these equations and was never 
explicitly clarified in the literature before, to the best of our knowledge.
Finally, numerical experiments demonstrate the favourable behaviour of the proposed numerical integrators. 
\end{abstract}

\section{Introduction}\label{sect-intro}
Since its introduction in the seminal paper \cite{hs91}, 
the Hunter--Saxton equation (HS equation below) 
\begin{align}\label{eqHSr-2}
u_{xxt}+2u_xu_{xx}+uu_{xxx} = 0, \quad x\in \mathbb{R} , \ t \geq 0
\end{align}
where $u:=u(x,t)$, has been attracting much attention.  
This is mainly due to its rich mathematical structures: 
the Hunter--Saxton equation is integrable; it is bihamiltonian; 
it possesses a Lax pair; it does not have global smooth solutions 
but enjoy two distinct classes of 
global weak solutions (conservative and dissipative);  
it can be seen as the geodesic equation of a right-invariant metric 
on a certain quotient space; etc. \cite{hz95,hz95b,bc05,l08} 
and references therein.  
Furthermore, the Hunter--Saxton equation arises 
as a model for the propagation of weakly nonlinear orientation waves 
in a nematic liquid crystal 
\cite{hs91} and it can be seen as the high frequency limit  
of another well known and well studied equation, namely 
the Camassa--Holm equation. More about this last equation can 
be found, for example, in the work \cite{raynaud}, 
the recent review \cite{hi10b}, 
and references therein.

Furthermore, the following extensions of the Hunter--Saxton equation enjoy 
intense ongoing research activities:  
the modified Hunter--Saxton equation (mHS equation below), introduced in \cite{l08b}, 
\begin{equation}\label{eqmHSr-2}
u_{xxt}+2u_xu_{xx}+uu_{xxx} -2\omega u_x= 0,
\end{equation}
where $\omega>0$, and the two-component 
Hunter--Saxton system (2HS below), introduced in \cite{w09}, 
\begin{align}
\begin{array}{l}
u_{xxt}+2u_xu_{xx}+uu_{xxx}-\kappa\rho\rho_x=0, \\
\rho_t+(u\rho)_x=0, 
\end{array}
\label{eq:HSsyst}
\end{align}
where $\kappa\in\{-1,1\}$ and $\rho:=\rho(x,t)$. 
These two partial differential equations (PDEs) also possess many interesting properties. 
The modified Hunter--Saxton equation is a model for short capillary 
waves propagating under the action of gravity \cite{l09}. An interesting 
feature of this modified version of the original problem is that it admits 
(smooth as well as cusped) travelling waves. This is not the 
case for the original problem \eqref{eqHSr-2}. Moreover, 
this PDE is also bihamiltonian \cite{l09}. 
The two-component generalisation of the Hunter--Saxton equation 
is a particular case of the Gurevich--Zybin system which describes 
the dynamics in a model of non-dissipative dark matter, 
see \cite{p05} and also \cite{lw12}. As the original equation, this system 
is integrable; has a Lax pair; is bihamiltonian; it is also the high-frequency limit of 
the two-component Camassa--Holm equation; has peakon solutions; its flow is equivalent to 
the geodesic flow on a certain sphere; etc. \cite{w10,k11,ml12,ww12,l13} 
and references therein.

Despite the fact that the above nonlinear PDEs are 
relatively well understood in a more theoretical 
way, there are not much results on numerical discretisations of these problems.
In fact, although we still continuously
find a number of papers on theoretical aspects
every year as of writing this paper, 
we are only aware of the numerical schemes from \cite{hkr07,xs09,xs10},
the latest being proposed in 2010.
All the three schemes are only for the original Hunter--Saxton equation.
The work \cite{hkr07} proves 
convergence of some discrete finite difference schemes to dissipative solutions of the Hunter--Saxton equation on the half-line. The references 
\cite{xs09,xs10} analyse local discontinuous Galerkin methods for the Hunter--Saxton 
equation and in particular, using results from \cite{hkr07}, prove convergence 
of the discretisation scheme to the dissipative solutions. 

The reason for such a lack of numerical studies should be
attributed to the following points.
First, due to the mixed derivatives present in the above problems, 
standard spatial discretisations of these nonlinear partial differential 
equations become, in general, nontrivial.
Second, partially in connection with this,
the Hunter--Saxton-like equations
are essentially underdetermined, and some strong (sometimes exotic)
assumptions are necessary to make the solution unique,
which causes difficulties in implementing numerical schemes.
In fact, in the existing numerical studies described above, 
some additional boundary conditions are employed to determine
the numerical solution without any systematic justification.

In the present publication, 
we will clarify the above issue and investigate two typical domain settings in details: 
the half real line and the periodic cases. 
We provide a systematic
view on the necessary and associated
additional boundary conditions under which the solution is determined
uniquely.
This not only gives a clear explanation for the additional conditions
employed in the existing schemes, but also provides a basis
for the new schemes in the present paper.

Based on this,
the main goal of this article is to present novel geometric numerical integrators for 
the Hunter--Saxton equation and for its two generalisations. 
As seen above, these nonlinear 
PDEs have important applications in physical sciences. 
One class of the proposed numerical schemes is based on new multi-symplectic formulations 
of the problems and are thus specially designed to preserve the multi-symplectic structures  
of the original equations. In addition, it was observed in the 
literature that multi-symplectic integrators have excellent potential 
for capturing long time dynamics of PDEs. 
Furthermore, the multi-symplectic schemes for the Hunter--Saxton like equations presented in this article are explicit integrators. 
The other class of schemes is based on Hamiltonian-preserving discretisations 
and are thus energy-preserving by construction. Most of these numerical schemes are implicit. 
A convergence analysis of the proposed numerical schemes 
will be reported elsewhere. 

We also illustrate the validity of the proposed schemes in some
numerical examples. 
More specifically, in this paper we focus on a typical exact solution
of the Hunter--Saxton equation that lacks smoothness,
and also several travelling waves for the other equations.
For both cases geometric numerical integrators are preferable;
for the former, geometric numerical integrators produce stable solutions without 
some stabilization such as upwinding
(we will show a comparison in
Section 3.1), and for the latter geometric numerical integrators have an advantage
in terms of long time behaviour.

The rest of the paper is organised as follows. 
We discuss the treatment of boundary conditions in the continuous setting with emphasis on
how we can determine solutions in
the Hunter--Saxton-like equations uniquely
in Section~\ref{sect-PI}. Then we 
present multi-symplectic formulations and numerical schemes 
for the considered class of PDEs in Section~\ref{sec-HShl}. The presentation 
of Hamiltonian-preserving numerical integrators is done in Section~\ref{sec-EP}. 
Finally, we draw some conclusions in Section~\ref{sec-conc}.

We will make use of the following notation. 
Spatial indices are denoted by $n$ and temporal ones by $i$. 
Numerical solutions are denoted by $u^n \simeq u( x_n, \cdot)$, 
$u^i \simeq u(\cdot, t_i )$
or $u^{n,i} \simeq u(x_n, t_i)$
on a uniform rectangular grid. We set $\Delta x=x_{n+1}-x_n, n\in\mathbb{Z}$, and
$\Delta t=t_{i+1}-t_i$, $i\geq 0$ . When we consider the domain $[-L,L]$, 
we set $x_0=-L, x_N = L$. When we consider periodic boundary conditions 
on the domain $[0,L]$, we set $x_0=0, x_N = L$, and assume $u^{n+N} = u^n $ as usual. 
We often write a solution 
as a vector $\bm u^i = (\ldots, u^{0,i},u^{1,i},u^{2,i},\dots )$ 
and use the abbreviation $\bm u^{i+1/2} = (\bm u^{i+1} + \bm u^i)/2$.
We define the forward and backward differences in time 
\begin{equation*}
\delta_t^+u^{n,i}=\frac{u^{n,i+1}-u^{n,i}}{\Delta t}
\qquad\hbox{and}\qquad 
\delta_t^-u^{n,i}=\frac{u^{n,i}-u^{n,i-1}}{\Delta t},
\end{equation*}
and similarly for differences in space. 
Also, we shall need the centered differences 
$\d_t=\hl(\delta_t^++\delta_t^-)$ 
and $\d_x=\hl(\delta_x^++\delta_x^-)$.
We also define the operator $\tilde{\delta}_x^2$ by 
\begin{equation*}
\tilde{\delta}_x^2 u^{n,i} = \dfrac{u^{n+1,i}-2u^{n,i}+u^{n-1,i}}{\Delta x^2}.
\end{equation*} 
Note that $\delta _x^2 \neq \tilde{\delta}_x^2$.
The trapezoidal rule is used as a discretisation of integrals:
\begin{equation*}
\dsum f^n \Delta x = \bra{\dfrac{1}{2}f^0 + \sum_{n=1}^{N-1}f^n + \dfrac{1}{2}f^N}\Delta x .
\end{equation*}
Finally, we will need the summation-by-parts formula
\begin{equation*}
\dsum f^n (\delta _x^+ g^n)\Delta x + \dsum (\delta _x^- f^n) g^n \Delta x 
=\left[
\frac{ f^n g^{n+1} + f^{n-1} g^n }{2}
\right]_0^N
\end{equation*}
which is frequently used for the analyses of Hamiltonian-preserving schemes. 

\section{On difficulties in solving numerically Hunter--Saxton-like equations: the treatment of boundary conditions}\label{sect-PI}

As mentioned in Section~\ref{sect-intro}, many difficulties arise in solving Hunter--Saxton-like (HS-like) equations numerically. 
Briefly speaking, solutions of HS-like equations are not unique 
due to the operator  $\px^2$ in front of $u_t$,
and thus even in the continuous setting it is a hard task
to determine how we should ``choose'' the solution.
Furthermore, even if it is done in the continuous case,
different difficulties arise in the discrete setting
where we are forced to impose additional discrete boundary conditions
(this will be illustrated and discussed in detail in the subsequent sections.)

In this section, we consider two typical domain settings,
the half real line case and the periodic case,
and describe the situation in detail in the continuous case.

In preparation for the following discussion,
let us first precisely
consider the meaning of ``underdetermined.''
Integrating both sides of the HS equation \eqref{eqHSr-2} twice
on some spatial domain, say $[-L, L]$ in view of numerical computation,
we obtain 
\begin{align}
u_{xxt} + 2u_xu_{xx} + uu_{xxx}=0, \label{a:r-2}\\
(u_t + uu_x)_x -\dfrac{1}{2}u_x^2 = a(t),  \label{a:r-1}\\
u_t + uu_x  - \px ^{-1} \left( \dfrac{1}{2}u_x^2 + a(t)\right) = h(t), \label{a:r0}
\end{align}
where $\px^{-1}(\cdot):=\int_{-L}^x (\cdot) \, \mathrm dx$, and $a(t)$ and $h(t)$ are 
the integral constants (independent of $x$). 
We refer to \eqref{a:r-2}, \eqref{a:r-1} and \eqref{a:r0} as ``rank $-2$'',
``rank $-1$'' and ``rank $0$'' equations, respectively
(``rank'' denotes $-1$ times the order of spatial differentiation of $u_t$).

Consider the rank $0$ equation \eqref{a:r0}.
It is reasonable to assume that the solution is unique for given $a(t)$ and $h(t)$ 
under appropriate boundary conditions.
If the constants and boundary conditions are given in advance,
one can discretise the rank $0$ equation.
However, if one considers the discretisation based on the lower rank equations 
a difficulty arises: solutions to the lower rank equations under the same boundary conditions as for the rank $0$ 
equation cannot be determined uniquely because of the lack of information on either $a(t)$ or $h(t)$. 
Hence, we must introduce some additional assumptions
to recover the lost information.
In what follows, we show how this can be done.

\subsection{The half line case}

The Hunter--Saxton equation \eqref{eqHSr-2} was originally
proposed as the rank $-1$ equation with $a(t)=0$
on the whole real line~\cite{hs91}.
There, with the method of characteristics,
it was shown that we need one boundary condition
``$u(x,t)\to 0$ as $x\to \infty$'' to determine 
a solution. 
Similarly, a solution can be determined under the boundary condition
$u(x,t)\to 0$ as $x\to -\infty$.
An exact (weak) solution can be constructed in closed form:
\begin{align} \label{hs:exsol}
u(x,t)=
\begin{cases}
0 & \text{if}\quad x\leq0\\
x/(0.5t+1) & \text{if}\quad 0<x<(0.5t+1)^2 \\
0.5t+1 & \text{if}\quad x\geq(0.5t+1)^2.
\end{cases}
\end{align}
See Figure~\ref{fig:ms} for a snapshot. 
In view of this solution, one typical setting of the Hunter--Saxton
equation is that we consider the equation on the half real
line $[0,\infty)$, and
we impose $u(0,t)=0$.
In this case, the unknown constant $h(t)$ can be easily found
 as $h(t)=0$ by considering the rank $0$ equation at $x=0$,
which is consistent with the claim that the solution can be
uniquely determined.
Many studies that follows~\cite{hs91}, including the numerical
studies~\cite{hkr07,xs09,xs10}, inherit this setting.

In the numerical computations, however, we have to ``cut'' 
the half line to a bounded domain, which we assume $[-L, L]$
without loss of generality.
Then it is natural to impose $u(-L,t)=0\ (t>0)$, which corresponds to
$u(x,t)\to 0$ as $x\to -\infty$.

Furthermore, when we consider finite difference discretisations, 
we discretise the equations uniformly in the spatial domain,
and thus the schemes use grid points outside the spatial domain.
Discrete boundary conditions for solving PDEs numerically are usually set 
such that they are consistent with the provided continuous boundary conditions.
This task is not straightforward for the HS equation \eqref{eqHSr-2},
since 
we generally  need boundary conditions at both the left and right boundaries,
while in the original half line setting
there is only one boundary condition $u(-L,t)=0$
(see also Remark~\ref{rem:onebc}).
Below we explain and justify the treatment in the existing numerical studies,
and how we can justify it.

First, in view of the exact solution above,
it is customary to impose $u_x(L,t)=0$, which is safe
for times $t$ in the time interval $[0,2(\sqrt{L}-1))$;
i.e., we understand that we take $L$ large enough so that 
the additional condition is safe for times we are interested in.

Depending on the finite difference scheme we employ, 
we often have to seek further additional conditions.
Following the existing studies and also
motivated by the discussion on the exact solution above,
let us focus on the solutions under the boundary conditions
$u(-L,t)=u_x(L,t)=0$ and with the constants $a(t)=h(t)=0$.
Note that $u_x(L,t)=0$ corresponds to $u_x(x,t) \to 0$ as $x\to\infty$
(for the theoretical studies under the condition $u_x(x,t) \to 0$ ($x\to\infty$), see, e.g.~Zhang--Zheng~\cite{zz00}).
We now assume that the solution of the rank $0$ equation \eqref{a:r0} in this problem setting
is unique
for each initial condition. 
Below we show that this unique solution 
automatically satisfies $u_x(-L,t) = u_{xx} (L, t) = 0$:

\begin{itemize}
\item $u_x (-L,t)=0$:\\
Let $v(\cdot) := u_x (-L,\cdot )$.
By evaluating \eqref{a:r-1} at $x=-L$, we obtain $\dot{v} + \dfrac{1}{2}v^2 = 0$.
In general, a solution of this ODE is of the form 
$v(t) = 2/ (t+c_1)$ with a constant $c_1$, but
as long as we consider an initial value satisfying $v(0) = 0$
(in Section~\ref{sec-HShl}, we considered such a case), we have $v(t)=0$ along the solution.
\item $u_{xx} (L,t) =0$:\\
Let $w(\cdot) := u_{xx} (L,\cdot )$.
By evaluating \eqref{a:r-1} at $x=L$, we obtain $uu_{xx} |_{x=L}=0$.
Here we exclude the case $u (L,t) = 0$, because this will give the solution $u(x,t)=0$.
This can be understood by \eqref{a:r0}, noting $\int _{-L}^L u_x^2\, \mathrm dx = 0$.
Therefore we immediately obtain $u_{xx}(L,t) = 0$.
\end{itemize}

These additional conditions were employed in the previous numerical studies~\cite{hkr07,xs09,xs10},
and are also used for geometric numerical integrators proposed in this paper.
The number of additional conditions depends on each scheme.

\begin{remark}\rm \label{rem:onebc}
If we carefully discretise the rank $-1$ equation,
it is possible to construct a scheme that happily works with
only one boundary condition $u(-L,t)=0$, which is consistent with
the continuous case.
As far as the present authors understand, it has never been pointed
out explicitly in the literature.
Such a scheme will be reported elsewhere since its exposition is outside the scope of this paper.
\end{remark}

\subsection{The periodic case}

First of all, we note an important fact that in the periodic case
the Hunter--Saxton-like equations
are essentially underdetermined in rank $-1$ (and accordingly rank $-2$),
since there is no way to ``add'' additional boundary conditions.
The only way to determine a solution is to provide the unknown constants
$a(t)$ and $h(t)$; see, for example, Yin~\cite{y04},
where the author considered the rank $-1$ HS with $a(t)=-(1/L)\int_0^L {u_x}^2/2\, {\rm d}x$ (which is necessary such that $u(x,t)$ is actually periodic)
and showed an unique existence theorem
under the assumption that $h(t)$ is explicitly given.
Thus one way to consider numerical computation is that
we seek a scheme that somehow incorporates the given information $h(t)$.

In the present paper, however, let us consider a different situation.
The mHS and 2HS equations have smooth travelling wave solutions,
which is in sharp contrast to the HS where any strong solutions
blow up in finite time~\cite{y04}.
Thus it makes sense to focus on such waves in the mHS and 2HS.

For the mHS equation, we focus on the smooth travelling waves 
of speed $c$, 
$u(x,t)=\varphi(x-ct)$, 
which are solutions to the differential 
equation \cite{l09}
\begin{align*}
(\varphi')^2=\frac{2\omega(M-\varphi)(\varphi-m)}{c-\varphi},
\end{align*}
where $M$, resp. $m$, is the maximum, resp. minimum, of the wave whenever $m<M<c$. 
The period of the wave is denoted by $L$, and the spatial domain is set to $[0,L]$.

Now we illustrate how we can ``recover'' the unknown constant $h(t)$
for such a travelling wave.
To see this, we first need to understand again why solutions  of the mHS 
equation \eqref{eqmHSr-2}, under periodic boundary conditions, are not unique. 
We proceed as for the HS equation. 
Integrating both sides of the rank $-2$ equation \eqref{eqmHSr-2}, we obtain the rank $-1$ equation
\begin{align}\label{eqmHSr-1}
(u_t + uu_x)_x -\dfrac{1}{2}u_x^2 - 2\omega u - a(t) = 0
\end{align}
with a constant $a(t)$. As long as we consider periodic boundary conditions, 
$a(t)$ has a unique expression of the form
\begin{align*}
a(t) = -\dfrac{1}{L} \int_0^L \left( \dfrac{1}{2}u_x^2 + 2\omega u\right)\, \mathrm dx, 
\end{align*}
because
\begin{align*}
0 = \int_0^L \left( (u_t + uu_x)_x -\dfrac{1}{2}u_x^2 - 2\omega u - a(t)\right) \mathrm dx = 
-\int_0^L \left( \dfrac{1}{2}u_x^2 + 2\omega u\right) \mathrm dx - L a(t).
\end{align*}
Furthermore, integrating \eqref{eqmHSr-1} once again, we obtain the rank $0$ equation
\begin{align*}
u_t + uu_x -\px^{-1} \left( \dfrac{1}{2}u_x^2 + 2\omega u + a(t)\right) = (u_t + uu_x) | _{x=0}
\end{align*}
where $\px^{-1}(\cdot) := \int _0^x (\cdot)\, \mathrm dx$.
By introducing $h(t) := (u_t + uu_x) | _{x=0}$, the above equation can be rewritten as
\begin{align}\label{eqmHSr0}
u_t + uu_x -\px^{-1} \left( \dfrac{1}{2}u_x^2 + 2\omega u + a(t)\right) = h(t).
\end{align}

Below we show an interesting fact that, as far as we are concerned with travelling waves 
under the periodic boundary condition, 
the constant $h(t)$, and accordingly
the uniqueness of the solution 
can be easily and automatically recovered by the concept of the pseudo-inverse 
of the differential (or difference) operator. Although pseudo-inverses are quite 
common and in fact often used in numerical analysis, it seems it has never 
been pointed out in the literature that pseudo-inverses can be
effectively utilised in such a way.

Since the concrete form of the mHS equation \eqref{eqmHSr0} 
is not essential for this discussion, we consider partial differential equations of the general form 
\begin{equation}\label{eq:formal}
u_t + f(u,u_x,u_{xx},\dots ) = h(t),
\end{equation}
as an equation of rank $0$ 
(below we often use the abbreviation $f(u) = f(u,u_x,u_{xx},\dots )$).
Assume that \eqref{eq:formal} has a periodic travelling wave solution
$u(x,t) = \varphi (x-ct)$ with $\varphi (x) = \varphi (x+L)$, 
where $L$ denotes the length of the period.
Substituting $\varphi$ into \eqref{eq:formal}, we obtain
\begin{equation*}
-c\varphi _x + f(\varphi, \varphi _x, \dots) = h(t).
\end{equation*}
For any $g$ satisfying $[g(\varphi)]_0^L=0$, one can re-express 
$h(t)$ as follows:
\begin{align*}
h(t) &=
\dfrac{\int_0^L g^\prime (\varphi)\, \dx}{\int_0^L g^\prime (\varphi)\, \dx} h(t)
+\dfrac{\int_0^L g^\prime (\varphi) c\varphi _x\, \dx}{\int_0^L g^\prime (\varphi)\, \dx} \\
&=
\dfrac{\int_0^L g^\prime (\varphi) (h(t) + c\varphi _x)\, \dx}{\int_0^L g^\prime (\varphi)\, \dx}
=
\dfrac{\int_0^L g^\prime (\varphi) f(\varphi )\, \dx}{\int_0^L g^\prime (\varphi)\, \dx}.
\end{align*}
As a special case if we select $g^\prime (\varphi) = 1$ (i.e., $g(\varphi)=\varphi$),
we have
\begin{equation} \label{eq:ht}
h(t) = \dfrac{\int_0^L f(\varphi) \dx}{L}.
\end{equation}
In short, if the equation \eqref{eq:formal} has a periodic travelling wave solution,
$h(t)$ should be expressed as \eqref{eq:ht}. 

Let us now consider the situation where
one does not know $h(t)$ explicitly and thus one is forced to work 
on equations of lower ranks, for example,
\begin{equation*}
u_{tx} + f_x(u,u_x,u_{xx},\dots ) = 0.
\end{equation*}
Let us introduce a pseudo-inverse operator of $\px$, 
denoted by $\px^\dagger$ (one can define a pseudo-inverse operator 
for a closed linear operator between two Hilbert spaces,
see~\cite{gr92,gr95} for example), 
and consider the following PDE 
\begin{equation}\label{pi:r-1}
u_t + \px^\dagger f_x (u,u_x,u_{xx},\dots)=0.
\end{equation}
\begin{prop}
The partial differential equation \eqref{pi:r-1} is equivalent to
\begin{equation} \label{eq:tw}
u_t + f(u,u_x,u_{xx},\dots ) = \dfrac{\int_0^L f(u) \, \dx}{L}.
\end{equation}
\end{prop}
This proposition indicates that the problem \eqref{pi:r-1} automatically catches 
travelling waves whenever they exist. 
Thus, for a proper numerical discretisation 
of travelling waves of the mHS equation, 
one should use the formulation \eqref{pi:r-1}.
\begin{proof}
Since $\mathrm{Ker}(\px) = \{ v \ | \ v = \mathrm{const.} \}$,
equation \eqref{pi:r-1} can be rewritten as
\begin{equation*}
u_t + f(u,u_x,u_{xx},\dots ) = k(t),
\end{equation*}
where $k(t)$ is a constant independent of $x$ and minimises $\| u_t \| _{L^2}$.
Since
\begin{align*}
\int _0^L u_t^2\, \dx &= \int _0^L (k(t)-f(u))^2\, \dx
= Lk(t)^2 - 2\left(  \int _0^L f(u)\, \dx \right) k(t) + \int _0^L f(u)^2\, \dx \\
&= L \left( k(t) - \dfrac{\int _0^L f(u) \dx}{L} \right) ^2
- \dfrac{\left( \int_0^L f(u)\, \dx \right) ^2}{L} + \int _0^L f(u)^2\, \dx,
\end{align*}
$k(t)$ should be of the form 
\begin{equation*}
k(t) = \dfrac{\int _0^L f(u)\, \dx}{L}.
\end{equation*}
This concludes the proof of the proposition. 
\qed
\end{proof}

A similar discussion for equations of rank $-2$ shows that the problem 
\begin{equation*}
u_t + (\px^2)^\dagger f_{xx} (u,u_x,u_{xx},\dots)=0
\end{equation*}
is equivalent to \eqref{eq:tw}.

In this paper, we also consider the periodic 2HS equation \eqref{eq:HSsyst} and its travelling waves. 
When $\kappa=1$, the smooth periodic travelling waves 
of speed $c$, $u(x,t)=\varphi(x-ct)$, are solutions 
to the differential equation \cite{ll09}
\begin{align} \label{2hs:diff}
(\varphi')^2&=\frac{b(Z-\varphi)(\varphi-z)}{c-\varphi}^2,
\end{align}
where $Z:=\max_{y\in\R}{\varphi(y)}$, resp. $z:=\min_{y\in\R}{\varphi(y)}$, is the maximum, 
resp. minimum, of the wave. Here, $b$ is a an additional positive parameter. Furthermore, 
\begin{align} \label{2hs:diff1}
\rho(x,t)=\psi(x-ct)=\frac{a}{c-\varphi},
\end{align}
where the parameter $a$ is determined by $b,z,Z$ and $c$, i.e. $a=\sqrt{b(c-z)(c-Z)}$.
Similar to the discussion about the periodic mHS equation seen above, one can show that  
a correct simulation of the travelling waves for the periodic 2HS equation should be based 
on the following differential equation
with the pseudo-inverse operator $(\px^2)^\dagger$:
\begin{eqnarray}
\begin{array}{l}
u_t+(\px^2)^\dagger (2u_xu_{xx}+uu_{xxx}-\kappa\rho\rho_x)=0, \\
\rho_t+(u\rho)_x=0.
\end{array}
\end{eqnarray}

\section{Multi-symplectic integrations of Hunter--Saxton-like equations}\label{sec-HShl}
We shall first present two new multi-symplectic formulations of the Hunter--Saxton equation, 
the corresponding explicit multi-symplectic schemes, as well as numerical experiments 
supporting our theoretical findings in Subsection~\ref{ms_hs}. 
Subsection~\ref{sec-mHSp} gives a similar program  for the multi-symplectic discretisation 
of the modified Hunter--Saxton equation. Using ideas from the above mentioned subsections, 
we shall present multi-symplectic schemes for the two-component Hunter--Saxton equation 
in Subsection~\ref{sec-2HS}.

The proposed multi-symplectic formulations of the HS-like equations follow the one presented 
in \cite{cor08} for the Camassa--Holm equation and in \cite{cmr14} 
for the two-component Camassa--Holm equation. 

For a detailed exposition on the concept of multi-symplectic partial differential equations, 
we refer the reader to, for example, the early references \cite{mps98,brid97,br01}. 

\subsection{Multi-symplectic discretisations of the Hunter--Saxton equation (half line case)}\label{ms_hs}
We consider the numerical discretisation of the Hunter--Saxton (HS) equation \eqref{eqHSr-2} 
in a computational domain $[-L,L] $ for a positive real number $L$ 
(correspondingly, we set $x_0=-L$ and $x_N=L$ with a natural number $N$).
In this subsection,  
we consider the boundary conditions given by 
$u(-L,t)=u_x (-L,t)=u_x(L,t)=0$
(see the previous section for the consistency of these conditions). 

\subsubsection{First multi-symplectic formulation and integrator for HS}\label{ms_hs1}
The multi-symplectic formulation 
\begin{equation}\label{eq:MS}
M\,z_t + K\,z_x = \nabla_z S(z)
\end{equation}
of \eqref{eqHSr-2} is obtained with $z=[u,\phi,w,v,\eta]^T$, 
the gradient of the scalar function
$S(z)=-w\,u-u\,\eta^2/2+\eta\,v$ and the 
two skew-symmetric matrices
\begin{equation*}
M=\begin{bmatrix}
0 & 0 & 0 & 0 & -\frac{1}{2} \\
0 & 0 & 0 & 0 & 0\\
0 & 0 & 0 & 0 & 0\\
0 & 0 & 0 & 0 & 0 \\
\frac{1}{2} & 0 & 0 & 0  & 0 
\end{bmatrix}
,\qquad
K=\begin{bmatrix}  
0 & 0 & 0 & -1 & 0 \\
0 & 0 & 1 & 0 & 0 \\
0 & -1 & 0 & 0 & 0 \\
1 & 0 & 0 & 0 & 0\\
0 & 0 & 0 & 0 & 0  
\end{bmatrix}. 
\end{equation*}  
For convenience, we also write this system 
componentwise
\begin{align*}
-\frac{1}{2}\eta_t-v_x
&=-w-\frac{1}{2}\eta^2,\\ 
w_x &= 0, \\ 
-\phi_x &= -u,\\
u_x &= \eta, \\ 
\frac{1}{2}u_t &=-u\eta+v. 
\end{align*}
A key observation, \cite{br01}, for the above multi-symplectic formulation of our problem 
is that the two skew-symmetric matrices $M$ and $K$ define symplectic structures on subspaces 
of $\R^5$
$$
\omega=\mathrm dz\wedge M\mathrm dz,\quad\quad \zeta=\mathrm dz\wedge K\mathrm dz,
$$
thus resulting in the following multi-symplectic conservation law 
\begin{equation}
\label{eq:msconslaw}
\partial_t\omega+\partial_x\zeta=0.
\end{equation}
This is a local property of our problem and we thus hope that multi-symplectic 
numerical schemes, as derived soon, 
will render well local properties of the HS equation \eqref{eqHSr-2}. 
More explicitly, we have for any solutions of \eqref{eq:MS}, the local conservation laws 
\begin{equation*}
\partial_t E(z)+\partial_x F(z)=0\qquad\hbox{and}\qquad
\partial_t I(z)+\partial_x G(z)=0,
\end{equation*}
with the density functions 
\begin{eqnarray*}
E(z)=S(z)-\frac{1}{2}z_x^TK^Tz\,,\qquad F(z)=\frac{1}{2}z_t^TK^Tz,\\
G(z)=S(z)-\frac{1}{2}z_t^TM^Tz\,,\qquad I(z)=\frac{1}{2}z_x^TM^Tz.
\end{eqnarray*}
Under the usual assumption on vanishing boundary terms for the functions
$F(z)$ and $G(z)$ one obtains the following global conserved quantities 
\begin{equation}\label{globinv}
\mathcal{E}(z)=\int_{-L}^L E(z)\,\mathrm dx\qquad\hbox{and}\qquad
\mathcal{I}(z)=\int_{-L}^L I(z)\,\mathrm dx. 
\end{equation}
For our choice of the skew-symmetric matrices $M$ and $K$, one thus 
obtains the density functions
\begin{align*}
E(z)=&S(z)+\dfrac{1}{2}z_x^TKz=\dfrac12(-wu-u(u_x)^2+(uv)_x),\\
F(z)=&-\dfrac{1}{2}z_t^TKz=\dfrac12u_tv-\dfrac12\phi_tw+\dfrac12w_t\phi-\dfrac12v_tu,\\
G(z)=&S(z)+\dfrac{1}{2}z_t^TMz=
-wu-\dfrac12\eta^2u+v\eta-\dfrac14\eta u_t+\dfrac14\eta_tu,\\
I(z)=&-\dfrac{1}{2}z_x^TMz=\dfrac{1}{4}
\bigl(u_x\eta-\eta_xu\bigr). 
\end{align*}
This will help us to derive the corresponding global invariants \eqref{globinv}. 

We first integrate the local conservation law $\partial_t I(z)+\partial_x G(z)=0$ 
over the spatial domain and obtain the Hamiltonian  
\begin{align}\label{ham1HS}
\mathcal{H}_1(u,u_x)=\frac12\int_{-L}^L u_x^2\,\mathrm dx.
\end{align}
Indeed for the above local conservation law, one has 
\begin{equation*}
0=\frac{1}{4}\frac{\mathrm{d}}{\mathrm{d}t}\int_{-L}^L
(u_x^2-uu_{xx})\,{\mathrm d}x
+\Bigl[-\frac34u_{xt}u+\frac14u_tu_x-u^2u_{xx}\Bigr]_{-L}^L.
\end{equation*}
Using one integration by parts and 
the vanishing boundary conditions for $u$ one thus gets 
\begin{equation*}
\frac12\frac{\mathrm{d}}{\mathrm{d}t}\int_{-L}^L
u_x^2\,{\mathrm d}x=0
\end{equation*}
which after integration gives the Hamiltonian \eqref{ham1HS}. 

Similarly, the second conservation law $\partial_t E(z)+\partial_x F(z)=0$ 
is linked to the Hamiltonian
\begin{align*}
\mathcal{H}_2(u,u_x)=\dfrac{1}{2}\int_{-L}^L uu_x^2 \mathrm dx 
\end{align*}
but, for our boundary conditions, this expression is not constant along solutions to our problem
\begin{align*}
\dfrac{\mathrm d}{\mathrm dt} \mathcal H_2 (u,u_x) = \left. \dfrac{1}{2} u_t ^2 \right|_{x=L}.
\end{align*}
Indeed, 
noting that $w\equiv0$ so that the multi-symplectic formulation \eqref{eq:MS}
is equivalent to the rank $-1$ equation with $a(t) = 0$,
one can simplify the density functions 
as follows
\begin{align*}
E(z)&=-\dfrac{uu_x^2}2+\dfrac12(uv)_x,\\
F(z)&=\dfrac12u_tv-\dfrac12v_tu.
\end{align*}
Integrating now the second conservation law, one gets 
\begin{align*}
0&=-\dfrac{\mathrm d}{\mathrm dt} \mathcal H_2 (u,u_x)+
\dfrac{\mathrm d}{\mathrm dt}\Bigl[\dfrac{uv}2\Bigr]_{-L}^L+
\Bigl[\dfrac{u_tv-v_tu}{2}\Bigr]_{-L}^L\\
&=-\dfrac{\mathrm d}{\mathrm dt} \mathcal H_2 (u,u_x)+
\left.\dfrac{u_tv+uv_t}{2} \right|_{x=L}+ \left.\dfrac{u_tv-v_tu}2 \right|_{x=L}
\end{align*}
which reduces to 
\begin{align*}
\dfrac{\mathrm d}{\mathrm dt} \mathcal H_2 (u,u_x) = \left. \dfrac{1}{2} u_t ^2 \right|_{x=L}.
\end{align*}

We now derive a numerical scheme based on 
the above multi-symplectic formulation of 
the Hunter--Saxton equation.

Following \cite{mr03}, one may obtain an integrator satisfying 
a discrete version of the multi-symplectic conservation law 
\eqref{eq:msconslaw} by introducing a 
splitting of the two matrices $M$ and $K$ 
in \eqref{eq:MS}, setting $M=M_++M_-$, 
$K=K_++K_-$ where $M_+^T=-M_-$ 
and $K_+^T=-K_-$. We will only consider the following matrices 
$M_+=\frac12M$ and $K_+=\frac12K$ for the splitting, 
keeping in mind that the above splitting 
of the matrices is not unique.  
The corresponding Euler box scheme reads
\begin{equation}\label{ms-euler}
M_+\delta_t^+z^{n,i}+M_-\delta_t^-z^{n,i}+K_+\delta_x^+z^{n,i}+
K_-\delta_x^-z^{n,i}=\nabla_zS(z^{n,i}),
\end{equation}
where $z^{n,i}\approx z(x_n,t_i)$.

The multi-symplecticity of the Euler box scheme is interpreted 
in the sense that, recall \eqref{eq:msconslaw}, 
\begin{equation*}
\delta_t^+\omega^{n,i}+\delta_x^+\zeta^{n,i}=0,
\end{equation*}
where $\omega^{n,i}=\text{d}z^{n,i-1}\wedge
M_+\text{d}z^{n,i}$ and
$\zeta^{n,i}=\text{d}z^{n-1,i}\wedge
K_+\text{d}z^{n,i}$.

With our choices for the matrices $M_+$ and $K_+$ and our multi-symplectic formulation 
\eqref{eq:MS} of the HS equation, the centered version 
of the Euler box scheme \eqref{ms-euler} 
reads (expressing the scheme only in the variable $u^{n,i}$) 
\begin{align*}
-\delta_x^2\delta_tu^{n,i}+\frac{1}{2}\delta_x\bigl((\delta_xu^{n,i})^2\bigr)-
\delta_x^2(u^{n,i}\delta_xu^{n,i})=0.
\end{align*}
Though the operator $\delta_x$ is not invertible, 
we will work on the following reformulation
\begin{align*}
-\delta_x\delta_tu^{n,i}+\frac{1}{2}(\delta_xu^{n,i})^2-
\delta_x(u^{n,i}\delta_xu^{n,i})=0,
\end{align*}
since it still retains the multi-symplecticity and is consistent with the rank $-1$ formulation.
Taking boundary conditions into account, we finally formulate the Euler box scheme as follows:
\begin{align}\label{EBHS}
v^{n,i}&:=\delta_xu^{n,i},\nonumber\\
\delta_tv^{n,i}&=\frac{1}{2}(v^{n,i})^2-
\delta_x(u^{n,i}v^{n,i})
\end{align}
for $n=1,\dots, N-1$, under the boundary conditions $u^{0,i}=0$, $u^{1,i} = u^{0,i}$,
$v^{0,i}=0$ and $v^{N,i}=0$.
They correspond to $u(-L,t)=0$, $u_x(-L,t)=0$, $u_x(-L,t)=0$ and $u_x(L,t)=0$,
respectively.
Let us point out 
that the numerical scheme \eqref{EBHS} is explicit, and since it preserves some geometry 
of the original PDE, it will perform well in terms of 
the evolutions of the Hamiltonians
as demonstrated in Subsection~\ref{ns_hsMS}.

\subsubsection{Second multi-symplectic formulation and integrator for HS}\label{ms_hs2}
As described in \cite{bc05}, see also \cite{cor08} 
in the context of the Camassa--Holm equation, 
in addition to \eqref{eqHSr-2} we 
can also consider the evolution equation satisfied by 
the energy density $\alpha:=u_x^2$. 
This permits to distinguish between two solutions, 
see \cite{bc05} for details. 
We thus obtain the following system of partial 
differential equations equivalent to the HS equation
\begin{subequations}
\begin{align*}
&u_t+uu_x+P_x=0,\\
&-P_{xx}=\frac12\alpha,\\
&\alpha_t+(u\alpha)_x=0. 
\end{align*}
\end{subequations}
The multi-symplectic formulation of the above system of partial 
differential equations is obtained setting 
$z=[u,\beta,w,\alpha,\phi,\gamma,P,r]$,
\begin{align*}
M=\begin{bmatrix}
0&-\frac12&0&0&0&0&0&0\\
\frac12&0&0&0&0&0&0&0\\
0&0&0&0&0&0&0&0\\
0&0&0&0&-\frac12&0&0&0\\
0&0&0&\frac12&0&0&0&0\\
0&0&0&0&0&0&0&0\\
0&0&0&0&0&0&0&0\\
0&0&0&0&0&0&0&0\\
\end{bmatrix},&&
K=\begin{bmatrix}
0&0&0&0&0&0&0&0\\
0&0&1&0&0&0&1&0\\
0&-1&0&0&0&0&0&0\\
0&0&0&0&0&0&0&0\\
0&0&0&0&0&1&0&0\\
0&0&0&0&-1&0&0&0\\
0&-1&0&0&0&0&0&-2\\
0&0&0&0&0&0&2&0\\
\end{bmatrix}
\end{align*}
and considering the scalar function 
\begin{equation*}
S(z)=-\gamma\,u+\frac{u^2\,\alpha}2-\alpha\,w+r^2. 
\end{equation*}
This is equivalent to the following system
\begin{align}\label{msHS2}
-\frac12\beta_t&=-\gamma+u\alpha,&&&
\frac12u_t+w_x+P_x&=0,\nonumber\\
-\beta_x&=-\alpha,&&&
-\frac12\phi_t&=-w+\frac{u^2}2,\nonumber\\
\frac12\alpha_t+\gamma_x&=0,&&&
-\phi_x&=-u,\nonumber\\
-\beta_x-2r_x&=0,&&&
2P_x&=2r. 
\end{align}

Using the above multi-symplectic formulation of 
our problem and similarly as before, 
the centered version of the Euler box scheme \eqref{ms-euler} 
for the second multi-symplectic formulation of 
the HS equation based on the choice of 
splitting matrices $M_+=\frac12M$ and $K_+=\frac12K$ reads
\begin{align}
\delta_t\alpha^{n,i}+\delta_x\bigl(u^{n,i}\alpha^{n,i}\bigr)&=0,\\
-\delta_{x}^2P^{n,i}&=\frac{1}{2}\alpha^{n,i}, \label{EBHS2}\\
\delta_tu^{n,i}+\frac{1}{2}\delta_x\bigl((u^{n,i})^2\bigr)+\delta_xP^{n,i}&=0.
\end{align}

We note that the associated global quantities \eqref{globinv} of the second multi-symplectic formulation 
are $\mathcal{H}_2$ and $\mathcal{H}_3$, see~\cite{hz94} for the definition of $\mathcal{H}_3$. 
But as will be explained soon, this multi-symplectic scheme \eqref{EBHS2} also offers a good behaviour for the 
evolution of $\mathcal{H}_1$. 

\subsubsection{Numerical simulations: Multi-symplectic schemes for HS}\label{ns_hsMS}
We now test our multi-symplectic integrators using the exact 
(nonsmooth) solution \eqref{hs:exsol} to the HS equation 
on the computational domain $[-6,6] $ (i.e., $L=6$) for times $0\leq t\leq T_{\text{end}}=0.5$.
Analytical values of the Hamiltonians are $\mathcal H_1 = 0.5$ and $\mathcal H_2 = t/8+1/4$.

\begin{figure}
\centering
\includegraphics*[height=4.3cm,keepaspectratio]{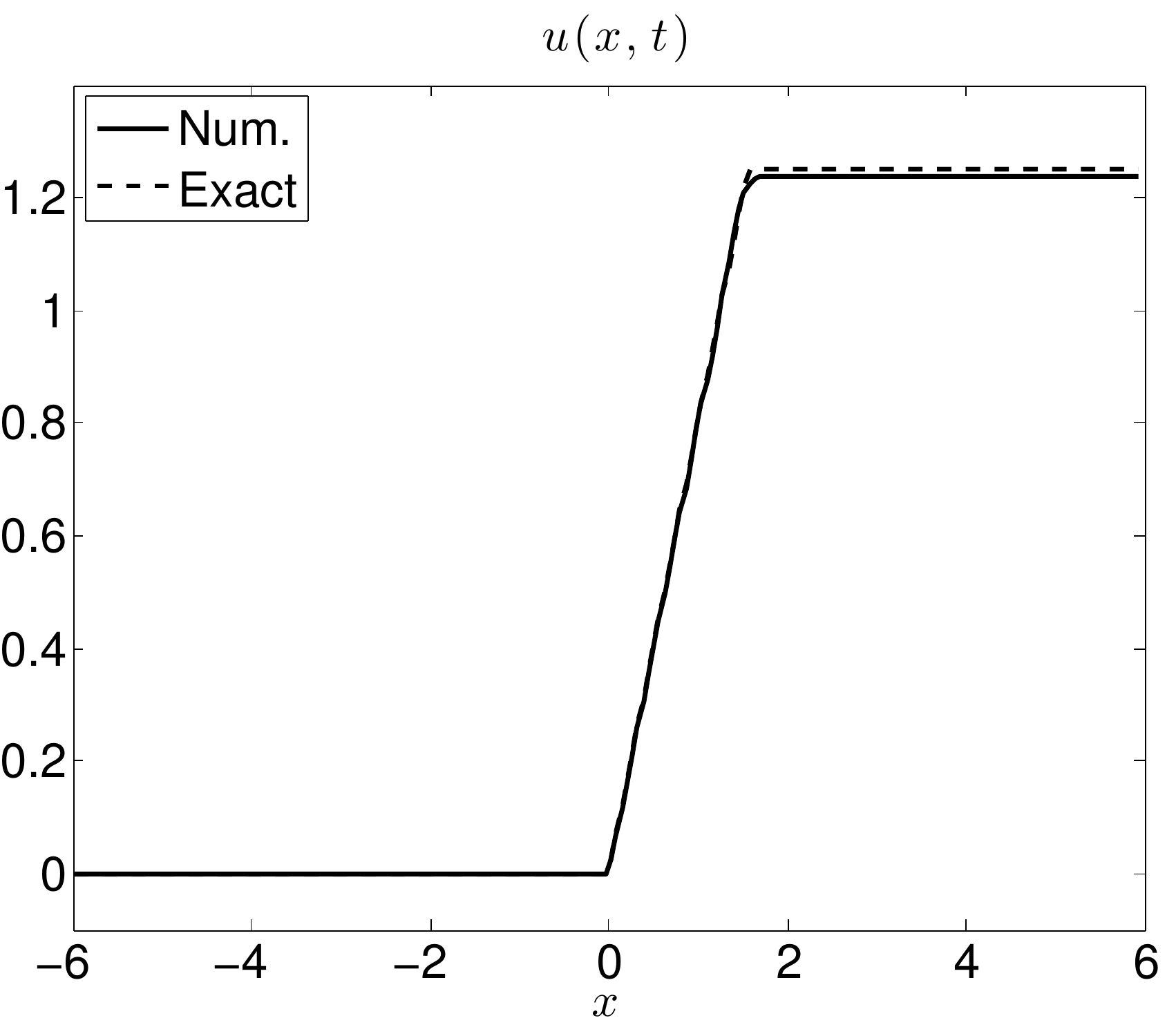}
\includegraphics*[height=4.3cm,keepaspectratio]{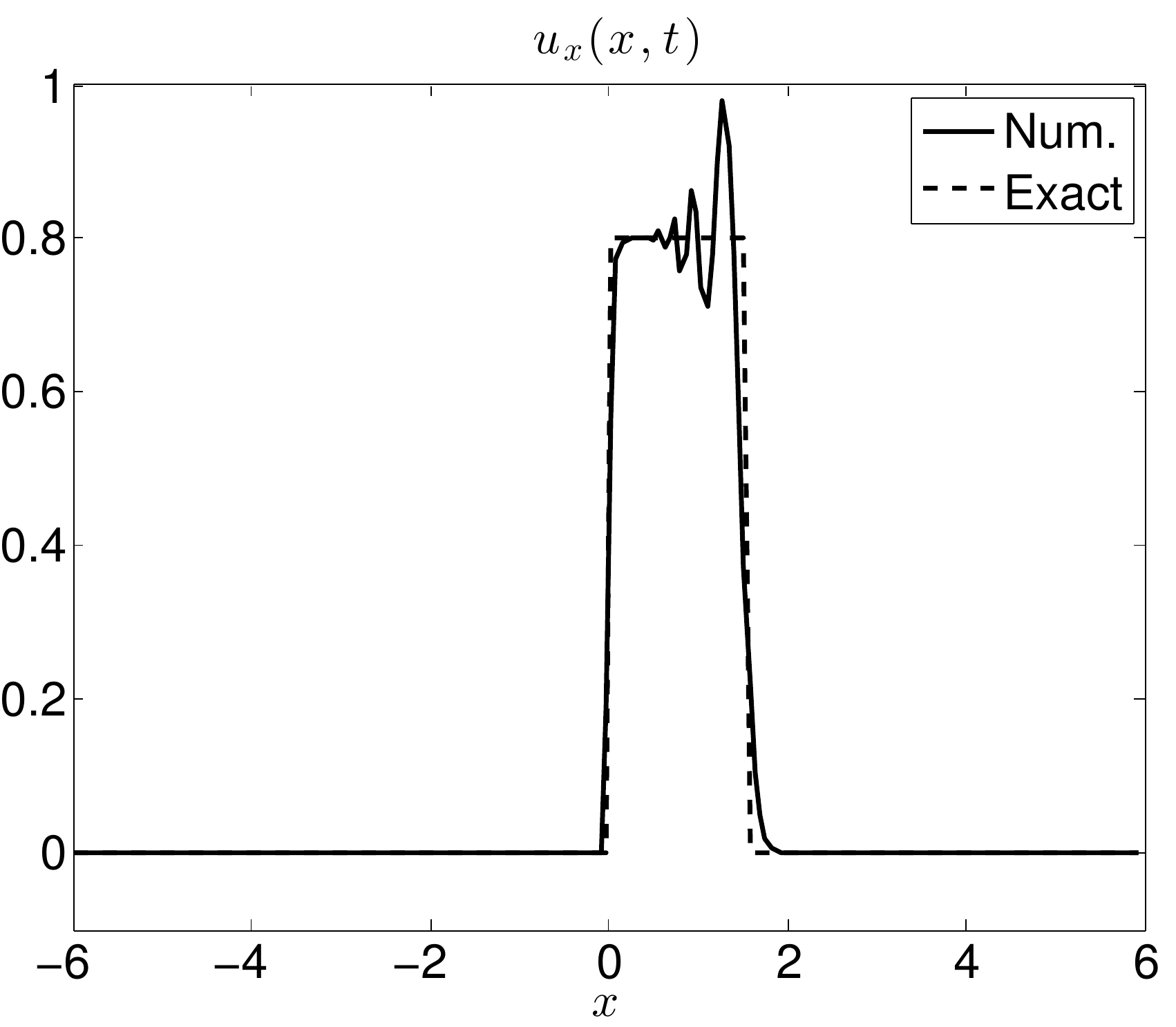}
\includegraphics*[height=4.3cm,keepaspectratio]{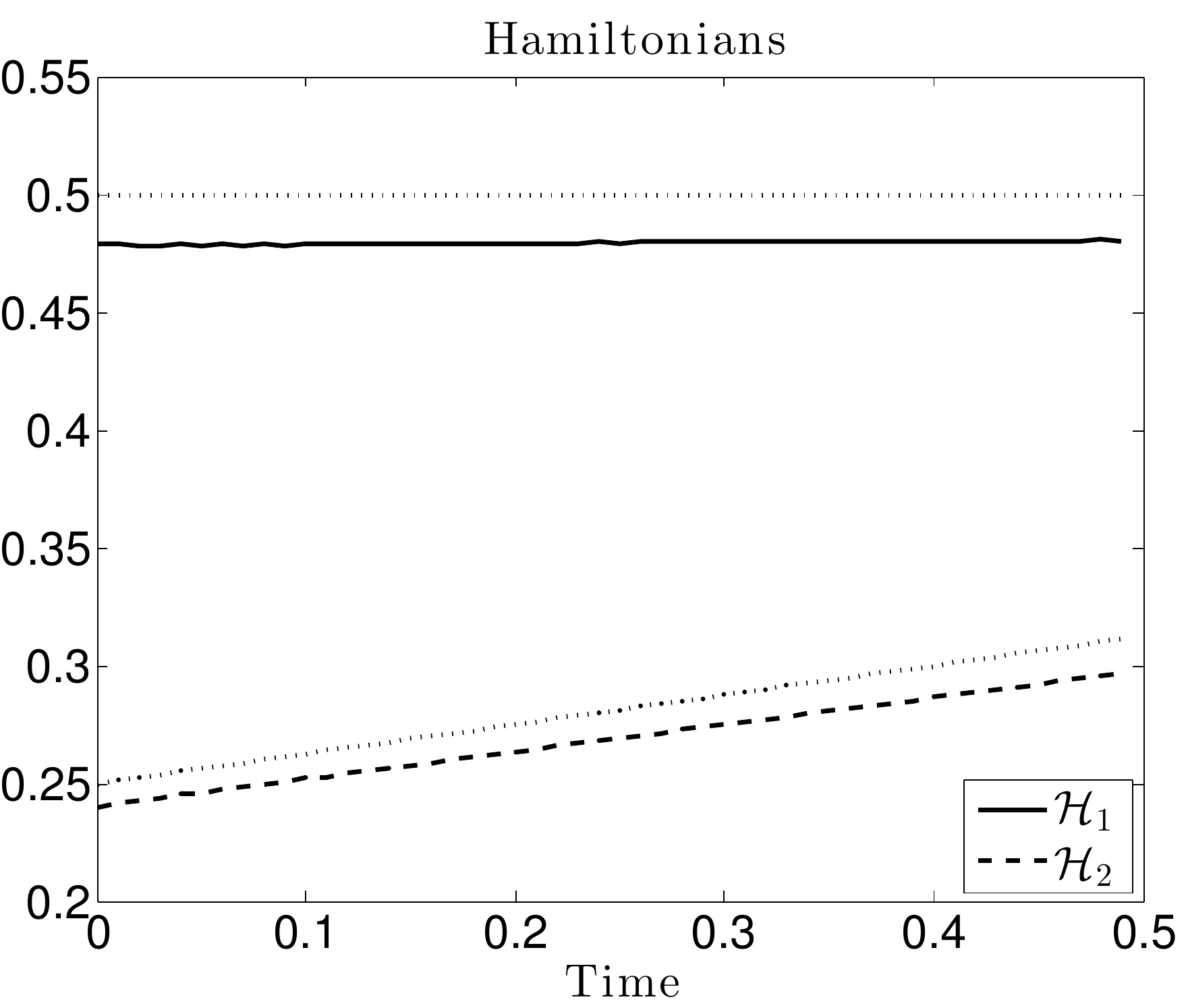}
\caption{Multi-symplectic scheme~\eqref{EBHS}: exact and numerical profiles of $u(x,t)$ and $u_x(x,t)$ at time $T_{\text{end}} = 0.5$
and computed Hamiltonians ($\Delta x = 12/201$ and $\Delta t=0.01$).}
\label{fig:ms}
\end{figure}

Figure~\ref{fig:ms} shows the numerical results obtained by 
the first multi-symplectic scheme \eqref{EBHS} with step sizes 
$\Delta x=12/201$ (i.e., $N=201$) and $\Delta t=0.01$.
The second multi-symplectic scheme \eqref{EBHS2} 
offers similar behaviours for both components of the solution and 
for both Hamiltonians as the first scheme \eqref{EBHS}. 
Results for this second numerical scheme are thus not displayed. 
As also observed for the numerical methods proposed in \cite{xs09},  
oscillations are present in $u_x(x,t)$, but 
the numerical approximation for $u(x,t)$, on compact intervals, is still correct. 
Moreover, one can observe the excellent evolutions of the Hamiltonians along the numerical solutions 
offered by the multi-symplectic scheme \eqref{EBHS}. 

\begin{figure}
\centering
\includegraphics*[height=4.3cm,keepaspectratio]{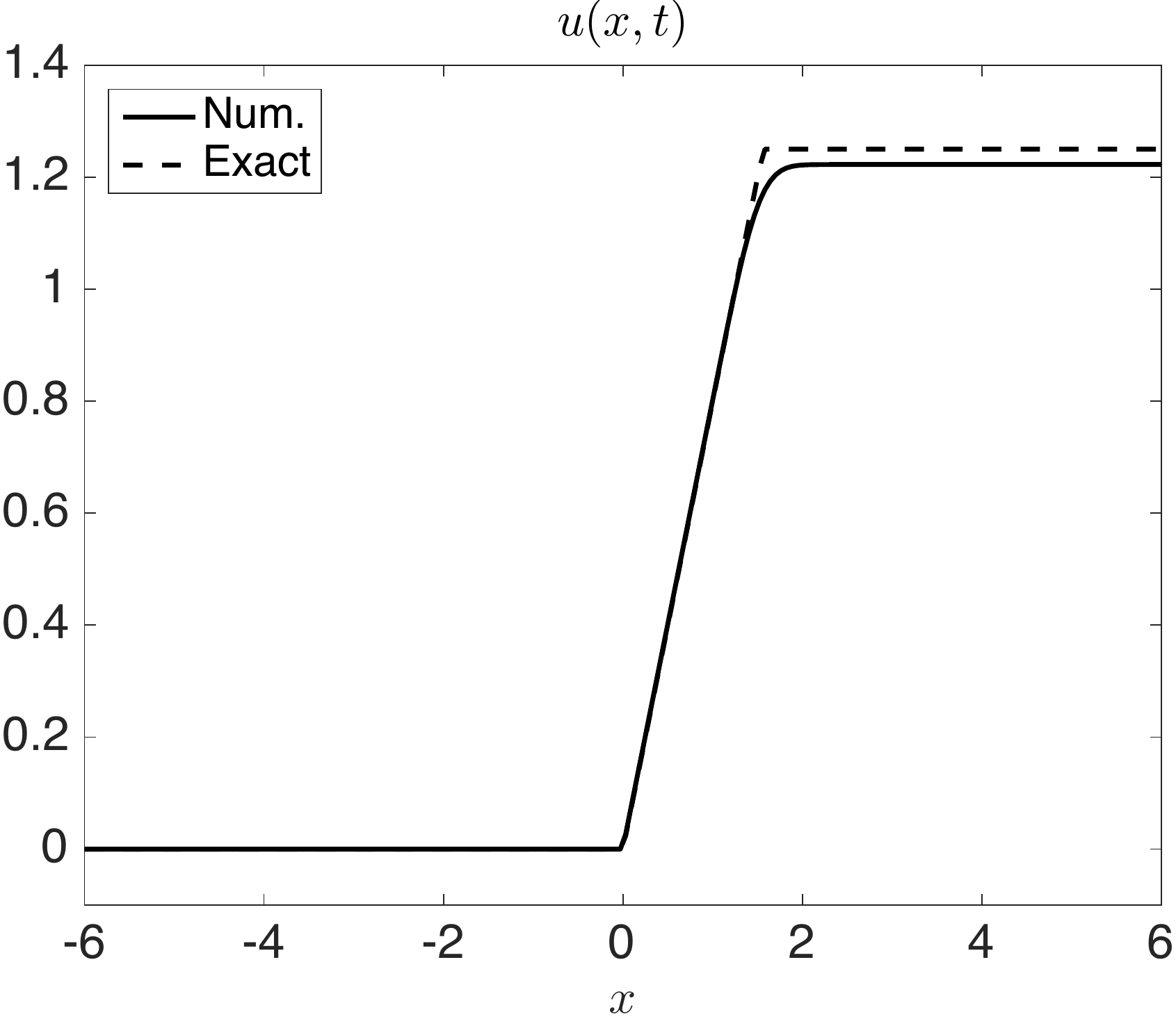}
\includegraphics*[height=4.3cm,keepaspectratio]{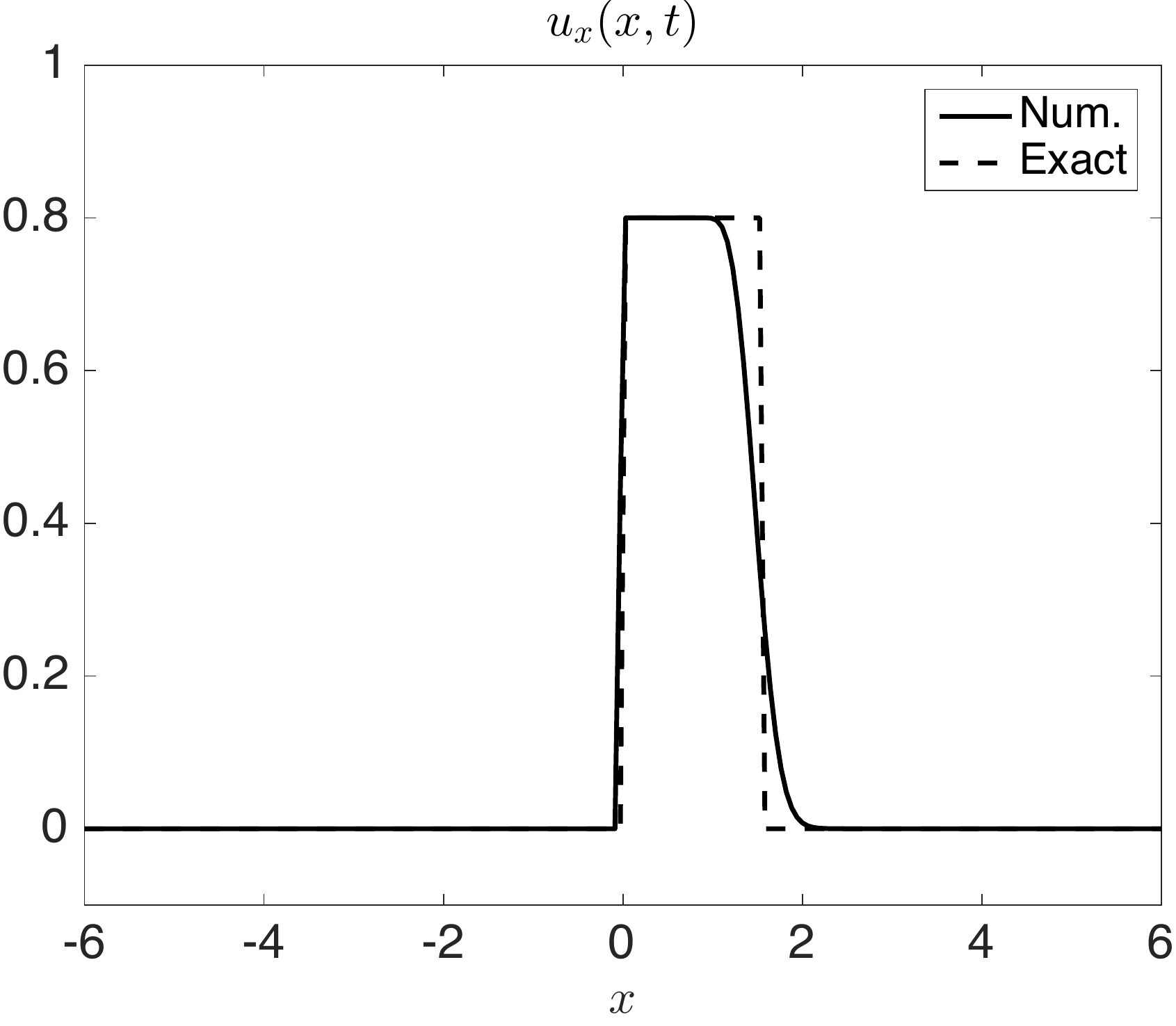}
\includegraphics*[height=4.3cm,keepaspectratio]{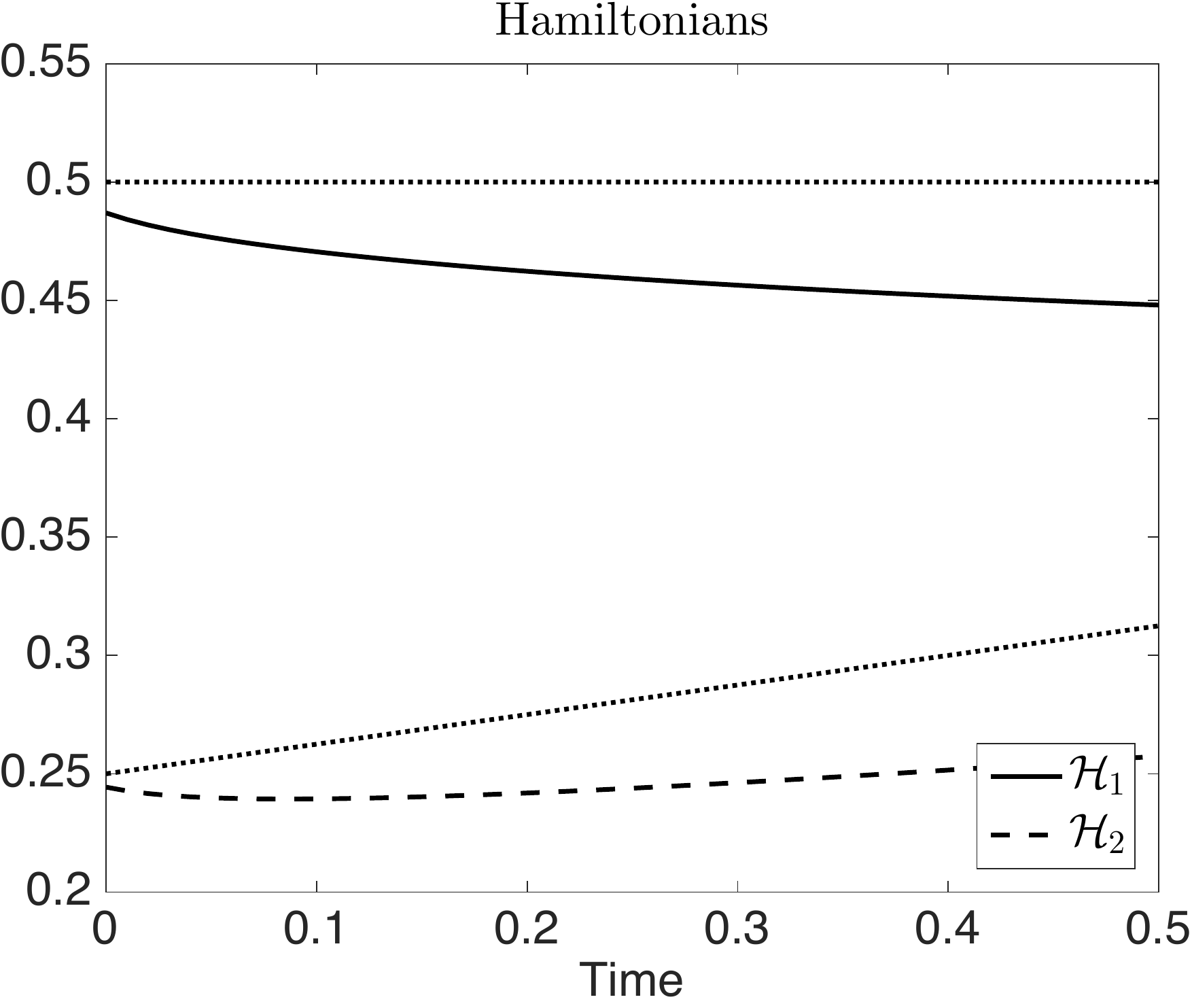}
\caption{Mid-point discretisation of the semi-discrete scheme by Holden et al. ((3.3) in~\cite{hkr07}): exact and numerical profiles of $u(x,t)$ and $u_x(x,t)$ at time $T_{\text{end}} = 0.5$
and computed Hamiltonians ($\Delta x = 12/201$ and $\Delta t=0.01$).}
\label{fig:hkr}
\end{figure}

Figure~\ref{fig:hkr} shows the numerical results obtained by the mid-point discretisation of the semi-discrete scheme proposed by Holden et al. ((3.3) in~\cite{hkr07}).
We compare our multi-symplectic scheme with this scheme.
Though no oscillation is observed for $u_x(x,t)$ in Figure~\ref{fig:hkr},
the multi-symplectic scheme \eqref{EBHS} seems preferable in terms of approximations of
$u(x,t)$ and the Hamiltonians.

\subsection{Multi-symplectic integrators for the modified Hunter--Saxton equation (periodic case)}\label{sec-mHSp}
In this section, we consider numerical integrators for the modified Hunter--Saxton (mHS) equation \eqref{eqmHSr-2}.
Note that for $\omega=0$, one obtains 
the Hunter--Saxton equation \eqref{eqHSr-2}.
An interesting feature of 
the mHS equation is that it admits 
(smooth as well as cusped) travelling waves when $\omega>0$.
This is not the case for the original problem \eqref{eqHSr-2}.
The aim of this subsection is thus to derive a multi-symplectic 
integrator for these travelling waves. 
Unfortunately, as seen in Section~\ref{sect-PI}, \eqref{eqmHSr-2} is an ``underdetermined'' 
PDE (in the sense that the problem has multiple solutions depending 
on the constant $h(t)$), which makes 
the construction of numerical schemes challenging. 
The results of Section~\ref{sect-PI} permit to overcome this difficulty 
and thus to derive a multi-symplectic 
scheme for the mHS equation. Finally, we present numerical experiments 
for the travelling waves of the mHS equation.

\subsubsection{Multi-symplectic formulation and integrator for mHS}
The following multi-symplectic formulation for the mHS equation follows from 
the one obtained for the HS equation in Subsection~\ref{ms_hs}. Indeed, the multi-symplectic formulation \eqref{eq:MS}
of the mHS equation is obtained with $z=[u,\phi,w,v,\eta]^T$, 
the gradient of the scalar function
$S(z)=-w\,u-u\,\eta^2/2+\eta\,v-\omega\,u^2$ and the 
two skew-symmetric matrices
\begin{equation*}
M=\begin{bmatrix}
0 & 0 & 0 & 0 & -\frac{1}{2} \\
0 & 0 & 0 & 0 & 0\\
0 & 0 & 0 & 0 & 0\\
0 & 0 & 0 & 0 & 0 \\
\frac{1}{2} & 0 & 0 & 0  & 0 
\end{bmatrix}
,\qquad
K=\begin{bmatrix}  
0 & 0 & 0 & -1 & 0 \\
0 & 0 & 1 & 0 & 0 \\
0 & -1 & 0 & 0 & 0 \\
1 & 0 & 0 & 0 & 0\\
0 & 0 & 0 & 0 & 0  
\end{bmatrix}. 
\end{equation*}  
For convenience, we also write this system 
componentwise
\begin{align*}
-\frac{1}{2}\eta_t-v_x
&=-w-\frac{1}{2}\eta^2-2\omega u,\\ 
w_x &= 0, \\ 
-\phi_x &= -u,\\
u_x &= \eta, \\ 
\frac{1}{2}u_t &=-u\eta+v. 
\end{align*}
As this was done in the previous section, we can integrate the first 
local conservation law $\partial_tI(z)+\partial_xG(z)=0$ over 
the spatial domain and obtain the Hamiltonian 
\begin{equation} \label{hamilmhs}
\mathcal{H}_1 (u,u_x) = \dfrac{1}{2} \int_0^L u_x^2\, \dx
\end{equation}
which is a conserved quantity for the mHS equation. 
The density functions for the second conservation law 
$\partial_tE(z)+\partial_xF(z)=0$ read
\begin{align*}
E(z)=-\dfrac{wu}2+\dfrac{(uv)_x}2-\omega u^2-\dfrac{uu_x^2}2,\\
F(z)=\dfrac12(u_tv-\phi_tw+w_t\phi-v_tu).
\end{align*}
Integrating this last conservation law gives 
\begin{align*}
0&=-\dfrac{\mathrm d}{\mathrm dt}\int_0^L(\dfrac{uu_x^2}2+\omega u^2)\,\mathrm dx
-\dfrac12\int_0^L(wu)_t\,\mathrm dx
+\dfrac12\dfrac{\mathrm d}{\mathrm dt}\int_0^L(uv)_x\,\mathrm dx\\
&\quad+\dfrac12\int_0^L(u_tv-\phi_tw+w_t\phi-v_tu)_x\,\mathrm dx.
\end{align*}
Observing that $w$ is constant in $x$ and $\eta,v$ are periodic, the above equation gives us 
\begin{align*}
\dfrac{\mathrm d}{\mathrm dt}\mathcal{H}_2:=\dfrac{\mathrm d}{\mathrm dt}\int_0^L(\dfrac{uu_x^2}2+\omega u^2)\,\mathrm dx=
-w\int_0^Lu_t\,\mathrm dx
\end{align*}
for the Hamiltonian
\begin{equation*}
\mathcal{H}_2 (u,u_x) = \dfrac{1}{2} \int_0^L (uu_x^2+2\omega u^2)\, \dx.
\end{equation*}
Note that the Hamiltonian $\mathcal H_2$ is a conserved quantity 
for travelling wave solutions but not for general solutions of the mHS equation. 

The explicit Euler box scheme \eqref{ms-euler} for the above multi-symplectic formulation reads
\begin{align*}
-\delta_x^2\delta_tu^{n,i}+\frac{1}{2}\delta_x\bigl((\delta_xu^{n,i})^2\bigr)-
\delta_x^2(u^{n,i}\delta_xu^{n,i})+2\omega\delta_xu^{n,i}=0.
\end{align*}
As seen in Section~\ref{sect-PI}, in order to select the right travelling waves, we have to consider 
the pseudo-inverse operator of $\delta _x^2$, denoted by $(\delta _x^2)^\dagger$. 
We then obtain the following numerical scheme  
\begin{align}\label{pEBmodHS}
-\delta_tu^{n,i}+(\delta _x^2)^\dagger \left( \frac{1}{2}\delta_x\bigl((\delta_xu^{n,i})^2\bigr)-
\delta_x^2(u^{n,i}\delta_xu^{n,i})+2\omega\delta_xu^{n,i} \right)=0.
\end{align}

\subsubsection{Multi-symplectic simulations of travelling waves for mHS}\label{mssimmHS}
\begin{figure}%
\begin{center}
\includegraphics*[scale=0.3]{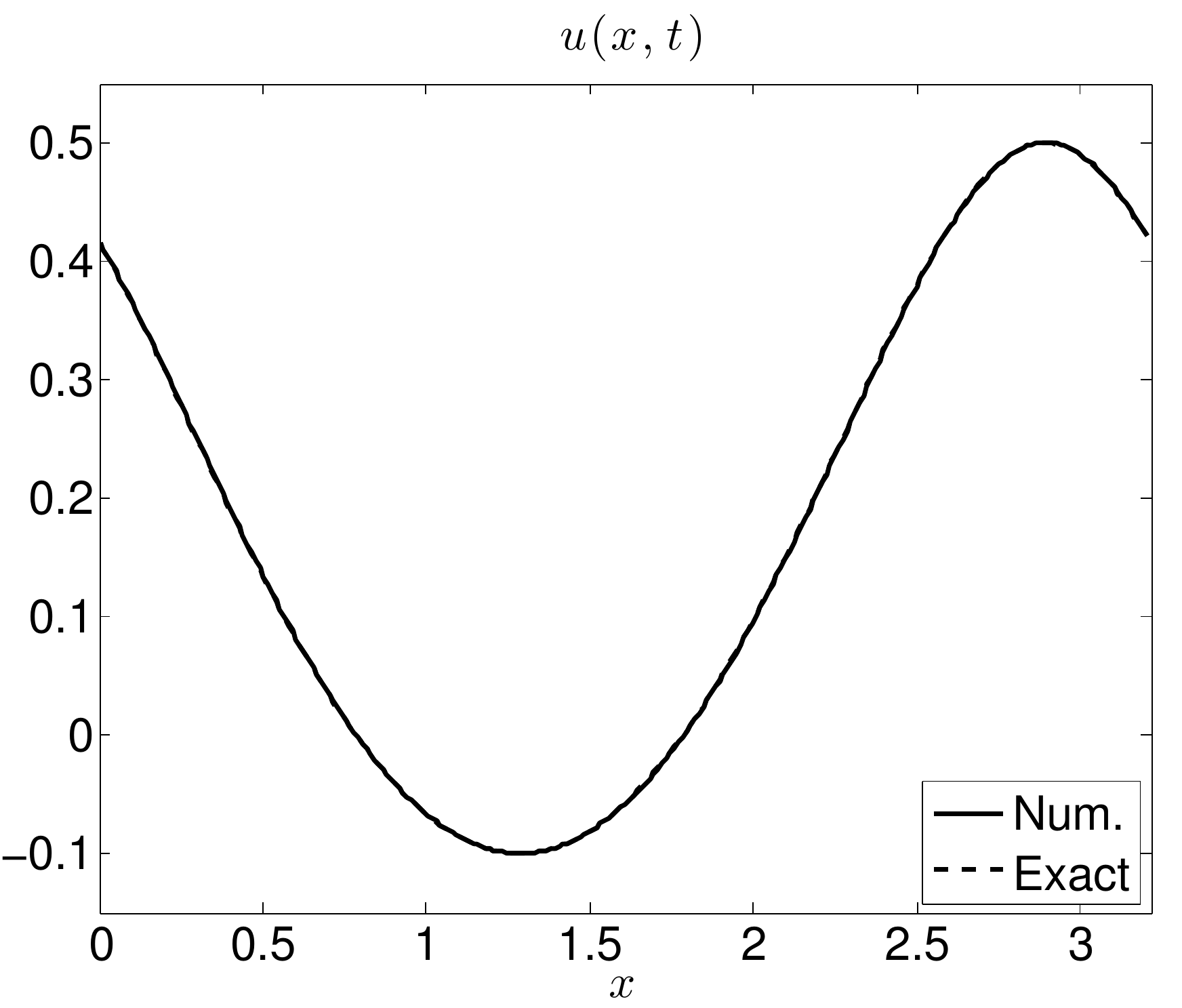}
\includegraphics*[scale=0.3]{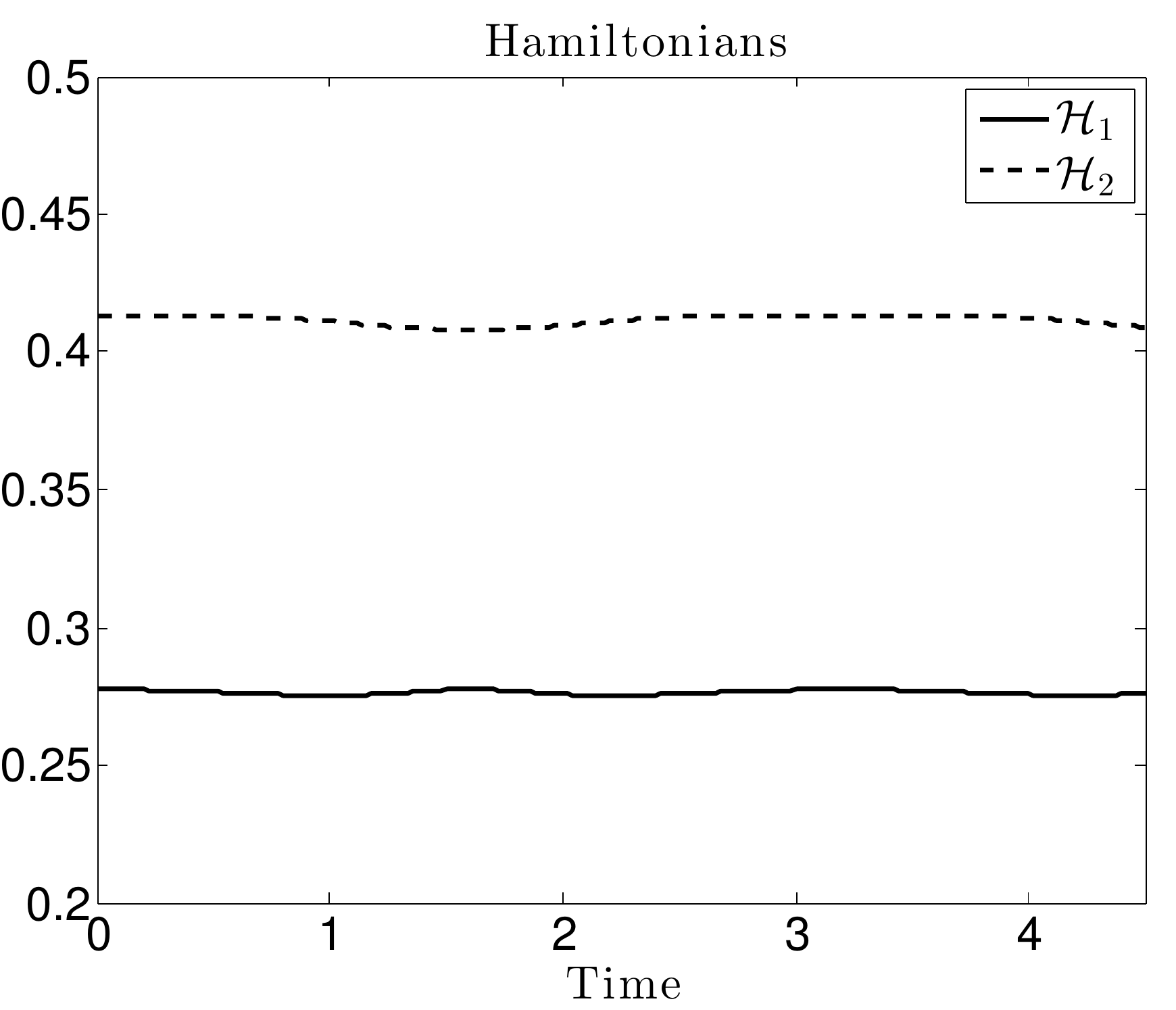}
\caption{Exact and numerical profiles of $u(x,t)$  
at time $T_{\text{end}}=3.5$ and computed Hamiltonians by the multi-symplectic scheme \eqref{pEBmodHS}. 
}
\label{fig:travmod}
\end{center}
\end{figure}

Figure~\ref{fig:travmod} displays the exact and 
numerical profiles of $u(x,t)$ at time $T_{\text{end}}=3.5$ 
and also the computed values of the Hamiltonians 
\begin{align*}
\mathcal{H}_1(u,u_x)=\frac12\int_0^L  u_{x}^2\,\mathrm{d}x,\quad
\mathcal{H}_2(u,u_x)=\frac12\int_0^L (uu_x^2+2\omega u^2)\,\mathrm{d}x
\end{align*} 
using the Euler box scheme \eqref{pEBmodHS} with (relative large) step sizes 
$\Delta t=0.02$ and $\Delta x=L_{\text{per}}/256$. 
We checked numerically that the solution 
of the above differential equation for $\varphi$ is periodic with period 
$L_{\text{per}}=3.2151\ldots$. The parameters for this 
simulation on the periodic computational domain $[0,L_{\text{per}}]$ are as follows: 
$\omega=1.5, m=-0.1, M=0.5, c=1$. One may notice that the numerical solution   
agrees very well with the exact one. Furthermore, good conservation 
properties of the numerical scheme are observed, even for such large step sizes, in the present figure. 

\subsection{Multi-symplectic discretisations of the two component Hunter--Saxton system (periodic case)}\label{sec-2HS}
In this section we consider the numerical discretisation of the two-component 
Hunter--Saxton (2HS) system \eqref{eq:HSsyst}   
introduced in~\cite{w09}. This system of partial differential equations also admits 
travelling wave solutions~\cite{ll09} so that 
one can use the pseudo-inverse as this was done in the previous subsection. 
We refer the reader to \cite{ll09} for an exposition of 
the Hamiltonian structures of the 2HS system.

\subsubsection{Multi-symplectic formulation and integrator for 2HS}
We can use the results from the preceding subsections 
to derive a multi-symplectic formulation for this system of equations too. 
Indeed, setting $z=[u,\phi,w,v,\eta,\rho,\gamma,\beta]$, and
using the following skew-symmetric matrices 
\begin{align}
M=\begin{bmatrix}
0 & 0 & 0 & 0 & -\frac12 & 0 & 0 & 0\\
0 & 0 & 0 & 0 & 0 & 0 & 0 & 0\\
0 & 0 & 0 & 0 & 0 & 0 & 0 & 0\\
0 & 0 & 0 & 0 & 0 & 0 & 0 & 0\\
\frac12 & 0 & 0 & 0 & 0 & 0 & 0 & 0\\
0 & 0 & 0 & 0 & 0 & 0 & 0 &\frac\kappa2\\
0 & 0 & 0 & 0 & 0 & 0 & 0 & 0\\
0 & 0 & 0 & 0 & 0 & -\frac\kappa2 & 0 & 0\\
\end{bmatrix},&&
K=\begin{bmatrix}
0 & 0 & 0 & -1 & 0 & 0 & 0 & 0\\
0 & 0 & 1 & 0 & 0 & 0 & 0 & 0\\
0 & -1 & 0 & 0 & 0 & 0 & 0 & 0\\
1 & 0 & 0 & 0 & 0 & 0 & 0 & 0\\
0 & 0 & 0 & 0 & 0 & 0 & 0 & 0\\
0 & 0 & 0 & 0 & 0 & 0 & 0 & 0\\
0 & 0 & 0 & 0 & 0 & 0 & 0 & \kappa\\
0 & 0 & 0 & 0 & 0 & 0 & -\kappa & 0\\
\end{bmatrix}
\label{mat:2hs}
\end{align}
and the gradient of the scalar function
\begin{equation*}
S(z)=-w\,u-u\,\eta^2/2+\eta\, v-\kappa\,u\,\rho^2/2+\kappa\,\gamma\,\rho
\end{equation*}
one obtains a multi-symplectic formulation \eqref{eq:MS} 
for the 2HS system \eqref{eq:HSsyst}. 
The above formulation reads componentwise:
\begin{align*}
-\frac{1}{2}\eta_t-v_x
&=-w-\frac{1}{2}\eta^2-\frac{\kappa}{2}\rho^2,\\ 
w_x &= 0, \\ 
-\phi_x &=-u,\\
u_x &= \eta, \\ 
\frac{1}{2}u_t &=-u\eta+v,\\
\frac{\kappa}{2}\beta_t &=\kappa\gamma-\kappa u\rho,\\
\kappa\beta_x &=\kappa\rho,\\
-\frac{\kappa}{2}\rho_t-\kappa\gamma_x &=0.
\end{align*}
The density functions used to compute the local conservation laws for 2HS are given by
\begin{align*}
E(z)&=-\dfrac12uw-\dfrac{uu_x^2}2+\dfrac12u_xv-\dfrac{\kappa}2u\rho^2+\dfrac{\kappa}2\gamma\rho
-\dfrac12w_x\phi+\dfrac12v_xu+\dfrac{\kappa}2\gamma_x\beta ,\\
F(z)&=\dfrac12\bigl(u_tv-\phi_tw+w_t\phi-uv_t-\kappa\gamma_t\beta+\kappa\beta_t\gamma\bigr) ,\\
I(z)&=\dfrac14\bigl(u_xu_x-u_{xx}u-\kappa\beta_{xx}\beta+\kappa\beta_x\beta_x\bigr) ,\\
G(z)&=u\eta^2-uv_x+\dfrac14u_t\eta-\dfrac14u\eta_t+\kappa\gamma\rho+\dfrac\kappa4\rho_t\beta-
\dfrac\kappa4\beta_t\rho.
\end{align*}
As before, using the fact that $w$ is constant in the variable $x$, 
$\eta$ and $v$ are periodic and integration by parts, the local conservation laws for 2HS 
give the following conserved quantities for travelling wave solutions
\begin{align}\label{hamil2HS}
\hspace{-1em}
\mathcal{H}_1(u,u_x,\rho)=\frac12\int(u_x^2+\kappa\rho^2)\,\mathrm{d}x, \quad
\mathcal{H}_2(u,u_x,\rho)=\frac12\int(\kappa u\rho^2+uu_x^2)\,\mathrm{d}x.
\end{align} 
We can now derive, as it was done previously, the centered version of 
the explicit Euler box scheme \eqref{ms-euler} for \eqref{eq:HSsyst} 
\begin{align}
\begin{array}{l}
-\delta_tu^{n,i}+(\delta_x^2)^\dagger\Bigl(\frac{1}{2}\delta_x\bigl((\delta_xu^{n,i})^2\bigr)-
\delta_x^2(u^{n,i}\delta_xu^{n,i})+\frac\kappa2\delta_x\bigl((\rho^{n,i})^2\bigr)\Bigr)=0, \\
\delta_t\rho^{n,i}+\delta_x(u^{n,i}\rho^{n,i})=0
\end{array}\label{EB2HS}
\end{align}
with the pseudo-inverse operator $(\delta_x^2)^\dagger$.

\begin{remark}
Similarly as above, one can also find a multi-symplectic formulation of the generalised periodic 
two-component Hunter--Saxton system \cite{ml12}
\begin{eqnarray}
\begin{array}{l}
u_{txx}+2\sigma u_xu_{xx}+\sigma uu_{xxx}-\kappa\rho\rho_x+Au_x=0, \\
\rho_t+(u\rho)_x=0. \end{array}
\label{eq:gHSsyst}
\end{eqnarray}
Here, in addition, $\sigma\in\R$ and $A\geq0$. 
Indeed, setting $z=[u,\phi,w,v,\eta,\rho,\gamma,\beta]$, 
using the skew-symmetric matrices \eqref{mat:2hs} 
and the gradient of the scalar function
\begin{equation*}
S(z)=-w\,u-\sigma\,u\,\eta^2/2+\eta\,v-\kappa\,u\,\rho^2/2+\kappa\,\gamma\,\rho+A\,u^2/2
\end{equation*}
one obtains a multi-symplectic formulation \eqref{eq:MS} 
for the generalised Hunter--Saxton system \eqref{eq:gHSsyst}. 
This would then permit to derive an Euler box scheme for these kind of equations too. 
\end{remark} 

\subsubsection{Multi-symplectic simulations of travelling waves for 2HS}\label{mssim2HS}
We now test the multi-symplectic Euler box scheme \eqref{EB2HS} on travelling 
wave solutions of \eqref{eq:HSsyst} with $\kappa=1$. 
The smooth periodic travelling waves 
of speed $c$, $u(x,t)=\varphi(x-ct)$, are solutions 
to the differential equation \eqref{2hs:diff}, see also \eqref{2hs:diff1}.

\begin{figure}%
\begin{center}
\includegraphics*[height=4.4cm,keepaspectratio]{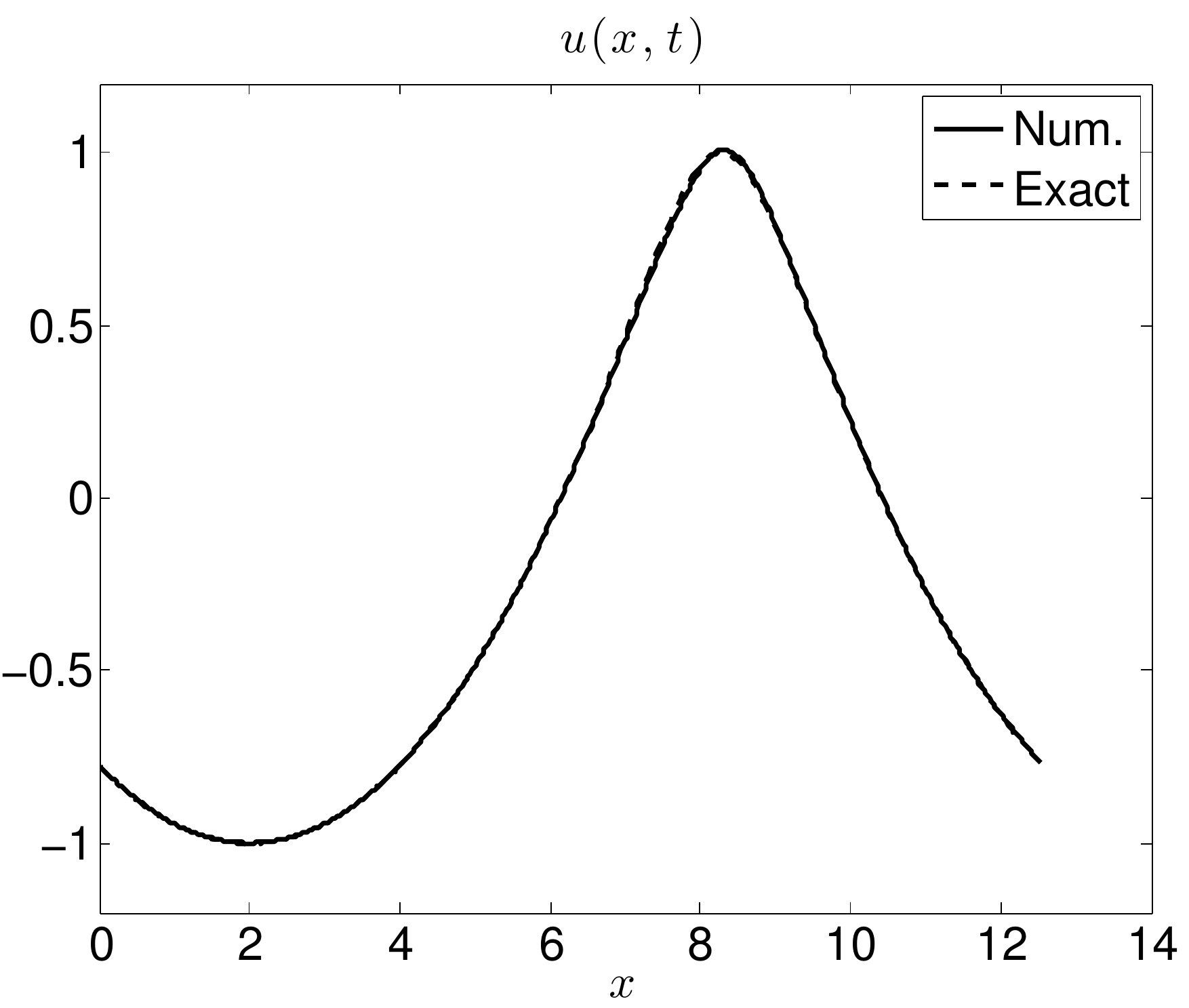}
\includegraphics*[height=4.4cm,keepaspectratio]{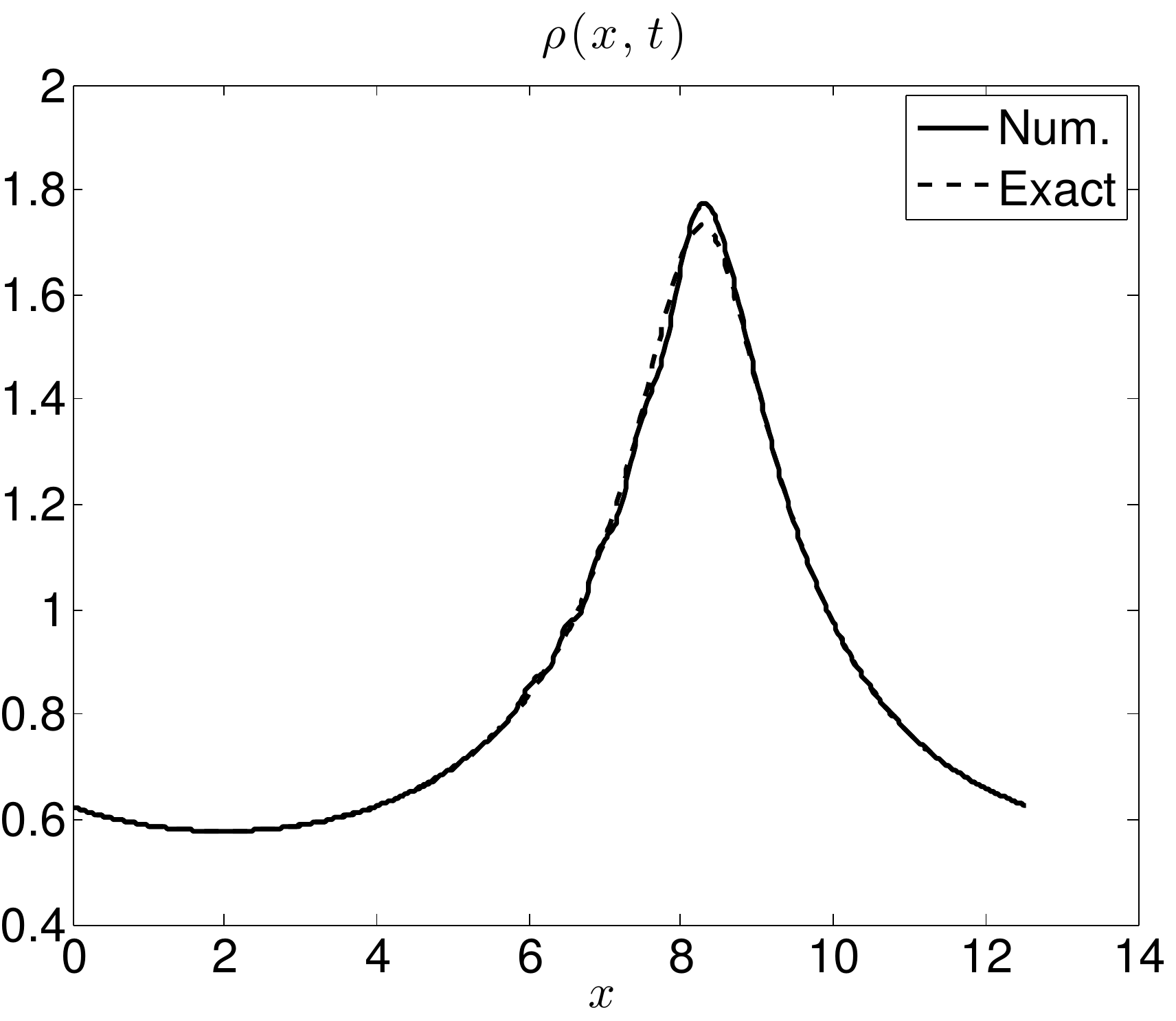}
\includegraphics*[height=4.4cm,keepaspectratio]{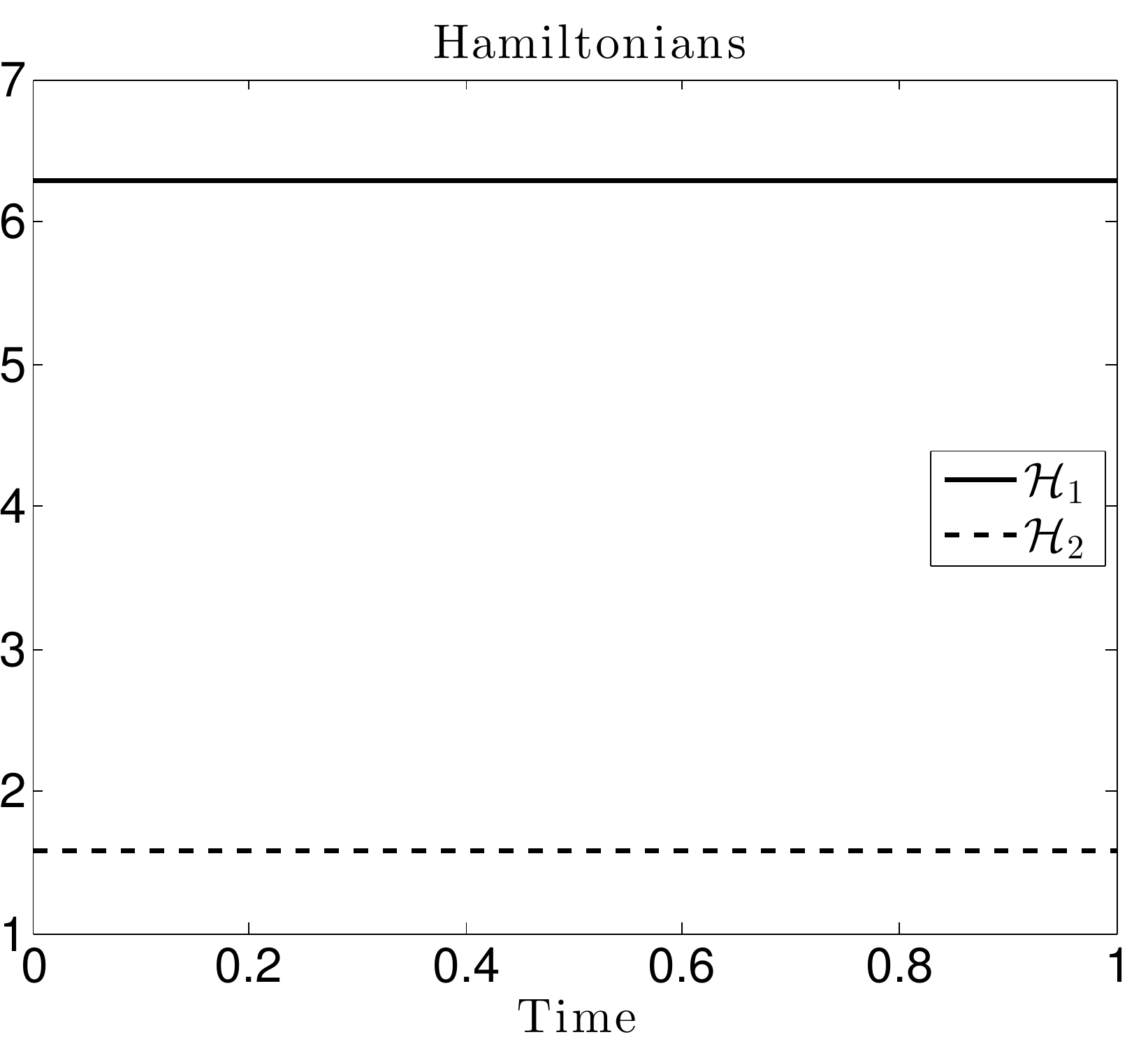}
\caption{Exact and numerical profiles of $u(x,t)$ and $\rho(x,t)$   
at time $T_{\text{end}}=1$ and computed Hamiltonians (right plot)
by the multi-symplectic scheme \eqref{EB2HS}. 
}
\label{fig:trav2}
\end{center}
\end{figure}

Figure~\ref{fig:trav2} displays the exact and 
numerical profiles of $u(x,t)$ and $\rho(x,t)$ at time $T_{\text{end}}=1$ 
and also the computed values of the Hamiltonians \eqref{hamil2HS} 
using the multi-symplectic Euler box scheme \eqref{EB2HS}  
with (large) step sizes 
$\Delta t=0.1$ and $\Delta x=L_{\text{per}}/512$.
The parameters for this 
simulation on the periodic computational domain $[0,L_{\text{per}}]$ are as follows: 
$z=-1, Z=1, b=1, c=2, a=\sqrt{3}$. We checked numerically that the solution of the differential 
equation \eqref{2hs:diff} for $\varphi(0)=z$ and $\varphi'(0)=0$ 
is periodic with period $L_{\text{per}}=12.5663\ldots$.  
One can observe, in this figure, that the numerical profile of the solution agrees very well with the exact profile. 
Furthermore, excellent conservation properties of the numerical solution are monitored. 

\section{Hamiltonian-preserving discretisations of Hunter--Saxton-like equations}\label{sec-EP}
As seen in the previous section, the class of HS-like equations, considered here, 
possesses invariants. It is thus of interest to derive invariant-preserving numerical methods. 
Furthermore, further numerical methods are useful for comparison purposes as exact solutions are 
generally not available for such problems. 

This section presents Hamiltonian-preserving schemes for the Hunter--Saxton equation, 
for the modified Hunter--Saxton equation as well as for the two-component Hunter--Saxton equation. 

Though the derivation of the proposed Hamiltonian-preserving schemes basically follows the lines of
\cite{fm11}, we need to pay attention on the treatment of boundary conditions.

\subsection{Hamiltonian-preserving integrators for HS (half line case)}
We first derive Hamiltonian-preserving numerical schemes for the HS equation based on the 
$\mathcal H_1$-variational formulation, resp. $\mathcal H_2$-variational formulation 
of the problem.
In this subsection, we assume the boundary conditions given by
$u(-L,t) = u_x (-L,t) = u_x (L,t) = u_{xx}(L,t) = 0$.

\subsubsection{An integrator based on the $\mathcal H_1$-variational formulation of HS}
\label{h1_hs}
In this subsection, we propose a finite difference scheme which preserves the Hamiltonian 
of the following variational formulation of the HS equation~\cite{hz94}:
\begin{align*}
u_{xxt} = (u_{xx} \px + \px u_{xx}) \px ^{-2} \fracdel{\mathcal H_1}{u},\quad
\mathcal H_1 (u,u_x) = \dfrac{1}{2} \int _{-L}^L u_x^2\, \mathrm dx,
\end{align*}
where $\px^{-1}(\cdot):=\int_{-L}^x (\cdot) \, \mathrm dx$.
As in the usual interpretation,
$(u_{xx} \px + \px u_{xx})$ operates on a function $f$ such that
$(u_{xx} \px + \px u_{xx})f = u_{xx}f_x + (u_{xx}f)_x$.
Since $\delta \mathcal H_1 / \delta u = -u_{xx}$, the above expression simplifies to 
\begin{align}\label{eqvfh1}
u_{xxt} = -(u_{xx}\px + \px u_{xx})u.
\end{align}
Under the boundary conditions $u(-L,t) = u_x (-L,t) = u_x (L,t) = u_{xx}(L,t) = 0$, we can prove the $\mathcal H_1$-preservation
based on \eqref{eqvfh1}.
In fact, using the integration-by-parts formula and \eqref{eqvfh1}, we have
\begin{align*}
\dfrac{\mathrm d}{\mathrm dt} \mathcal H_1 &=
\int _{-L}^L u_x u_{xt} \mathrm \, dx = - \int _{-L}^L uu_{xxt} \, \mathrm dx+ [uu_{xt}]_{-L}^L \\
&= \int _{-L}^L u \cdot (u_{xx}\px + \px u_{xx})u \, \mathrm dx = [u^2 u_{xx}]_{-L}^L = 0.
\end{align*}
Here, the boundary terms $[uu_{xt}]_{-L}^L$ and $[u^2 u_{xx}]_{-L}^L$
vanish due to the boundary conditions.

We now derive a finite difference scheme which preserves the Hamiltonian structure,
namely, an $\mathcal H_1$-preserving scheme. This derivation is based on 
the discrete variational derivative method~\cite{f99,fm11} (see also \cite{cgm12}). 
We focus on the spatial discretisation, since the idea of the temporal 
discretisation is essentially the same as the discrete gradient method 
(see, for example, \cite{g96,qm08} and references therein).

Firstly, let us define a discrete version of $\mathcal H_1$ by
\begin{equation*}
\mathcal H_{1, \rm d} (\bm u) :=
\dsum \dfrac{\Delta x}{2} \dfrac{(\delta_x^+ u^n)^2 + (\delta_x^- u^n)^2}{2}.
\end{equation*}
Based on the first Hamiltonian structure \eqref{eqvfh1},
we define the following semi-discrete scheme: 
\begin{equation}\label{h1sds}
\tilde{\delta}_x^2 \dot{u}^n =  - (\tilde{\delta}_x^2 u^n) (\delta_x u^n) - \delta_x (u^n (\tilde{\delta}_x^2 u^n))
\end{equation}
for $n=1,\dots,N$,
where the dot on $u$ stands for the differentiation with respect to time. We consider 
the numerical boundary conditions
$u^{0} = 0$, $u^{-1}=u^{1}$, $u^{N+1} = u^{N-1}$, $u^{N+2} = 2u^N - u^{N-2}$,
which correspond to $u(-L,t)=0$, $u_x(-L,t)=0$, $u_x(L,t)=0$, $u_{xx}(L,t)=0$, respectively.
The semi-discrete scheme \eqref{h1sds} inherits the $\mathcal H_1$-preservation.
To simplify the calculations of the proof, we extend the semi-discrete scheme to $n=0$.
Note that the scheme refers $u^{-2}$ when $n=0$, but we simply understand
$u^{-2}$ takes a value 
such that the scheme holds even for $n=0$. 
The $\mathcal H_1$-preservation is now checked as follows.
\begin{align*}
&\dfrac{\mathrm d}{\mathrm dt} \mathcal H_{1, \rm d} (\bm u(t))\\
&=
\dsum \dfrac{\Delta x}{2}
\left( (\delta_x^+ u^n)(\delta_x^+ \dot{u}^n) + (\delta_x^- u^n)(\delta_x^- \dot{u}^n) \right)\\
&=
-\dsum \Delta x u^n (\tilde{\delta}_x^2 \dot{u}^n)  \\
& \qquad +\frac{1}{4} \left[ (\delta_x^+ \dot{u}^n) u^{n+1} + (\delta_x^+ \dot{u}^{n-1}) u^{n}
 + u^{n} (\delta_x^- \dot{u}^{n+1}) + u^{n-1} (\delta_x^- \dot{u}^{n}) \right] _{0}^N \\
&=
\dsum \Delta x \left( u^n (\tilde{\delta}_x^2 u^n) (\delta_x u^n) + u^n \delta_x (u^n (\tilde{\delta}_x^2 u^n)) \right) \\
&= 
\dsum \Delta x \left( u^n (\tilde{\delta}_x^2 u^n) (\delta_x u^n) - (\delta _x u^n) u^n (\tilde{\delta}_x^2 u^n) \right) \\
&\phantom{=}
+\dfrac{1}{4} \left[
u^n \left( u^{n+1} (\tilde{\delta}_x^2 u^{n+1}) + u^{n-1} (\tilde{\delta}_x^2 u^{n-1}) + u^{n+1} (\tilde{\delta}_x^2 u^n) + u^{n-1} (\tilde{\delta}_x^2 u^n) \right)
\right] _{0}^N \\
&= 0.
\end{align*}

Next, we apply the discrete gradient method for the temporal discretisation in order to obtain
the following fully-discrete scheme:
\begin{equation}\label{h1fds}
\tilde{\delta}_x^2 \delta_t^+ u^{n,i} =  - (\tilde{\delta}_x^2 u^{n,i+1/2}) (\delta_x u^{n,i+1/2}) - \delta_x (u^{n,i+1/2} (\tilde{\delta}_x^2 u^{n,i+1/2}))
\end{equation}
for $n=1,\dots,N$,
still with the numerical boundary conditions
$u^{0,i} = 0$, $u^{{-1},i}=u^{1,i}$, $u^{N+1,i} = u^{N-1,i}$, $u^{N+2,i} = 2u^{N,i} - u^{N-2,i}$,
which correspond to $u(-L,t)=0$, $u_x(-L,t)=0$, $u_x(L,t)=0$, $u_{xx}(L,t)=0$, respectively,
for all $i$ and $u^{n,i+1/2}$ denotes an abbreviation for $(u^{n,i}+u^{n,i+1})/2$.
The scheme \eqref{h1fds} has the following conservation property by construction
\begin{equation*}
\mathcal H_{1,\rm d} (\bm{u}^i) = \mathcal H_{1,\rm d} (\bm{u}^0) \quad
{\rm for \ all} \ i. 
\end{equation*}
 
\begin{remark}
The resulting scheme strongly depends on the definition of the discrete version of $\mathcal H_1$.
Other definitions such as, for example, 
\begin{equation*}
\mathcal H_{1, \rm d} (\bm u) :=
\dsum \dfrac{\Delta x}{2} (\delta_x u^n)^2 
\end{equation*}
will lead to different $\mathcal H_1$-preserving schemes.
\end{remark}

\subsubsection{An integrator based on the $\mathcal H_2$-variational formulation of HS}
\label{h2_hs}
In this subsection, we propose a finite difference scheme 
which preserves the second Hamiltonian structure of the HS equation~\cite{hz94}:
\begin{align*}
u_{xt} = \fracdel{\mathcal H_2}{u},\quad \mathcal H_2 (u,u_x) = 
\dfrac{1}{2} \int _{-L} ^L uu_x^2\, \mathrm dx,
\end{align*}
or equivalently 
\begin{align*}
u_{xt} = -(uu_x)_x + \dfrac{u_x^2}{2}, \quad \text{or}\quad
u_t = -uu_x + \px ^{-1} \left( \dfrac{u_x^2}{2}\right) ,
\end{align*}
where $\px ^{-1}(\cdot) = \int _{-L}^x (\cdot)\, \mathrm dx$.
Note that $\mathcal H_2$ is not an invariant as 
already seen in Subsection~\ref{ms_hs}. 
This can also be confirmed using the above variational structure. 
In fact, we have
\begin{align}
\dfrac{\mathrm d}{\mathrm dt} \mathcal{H}_2
&=
\int _{-L}^L \left( \dfrac{1}{2}u_x^2 u_t + uu_x u_{xt} \right) \, \dx 
=
\int _{-L}^L \left( \dfrac{1}{2}u_x^2 - (uu_x)_x \right) u_t \, \dx + [uu_xu_t]_{-L}^L \nonumber \\
&=
\int _{-L}^L  u_{tx}  u_t \, \dx 
=
\left[ \dfrac{1}{2} u_t ^2 \right] _{-L}^L  
=
\left. \dfrac{1}{2}u_t ^2 \right| _{x=L}.
\label{H2_td}
\end{align}
Let us now define a discrete version of $\mathcal H_2$ by
\begin{equation*}
\mathcal H_{2, \rm d} (\bm u, \bm v) :=
\dsum \dfrac{\Delta x}{2} u^n (v^n)^2,
\end{equation*}
where $v^n:= \delta _x u^n$. 
Based on the second Hamiltonian structure,
we define the following semi-discrete scheme:
\begin{equation} \label{h2sds}
\dot{u}^n = - u^n v^n + \delta _x^{-1} \bra{\dfrac{1}{2}(v^n)^2}
\end{equation}
for $n=1,\dots,N$.
We assume the numerical boundary conditions $u^{0} = 0$ and $v^{N} = 0$,
which correspond to $u(-L,t)=0$ and $u_x(L,t)=0$, respectively,
and $\delta _x^{-1} \bra{\dfrac{1}{2}(v^n)^2}$ is defined by
\begin{align*}
\delta_x^{-1} \bra{\dfrac{1}{2}(v^n)^2} = 
\begin{cases}
0, & \text{if } n = 0,1, \\
 {\displaystyle
2 \Delta x \sum_{k=1}^m \frac{v_{2k-1}^2}{2}, }  & \text{if } n=2m \ (m\geq 1), \\
\displaystyle{ 2\Delta x \sum_{k=1}^m \frac{v_{2k}^2}{2},}  & \text{if } n=2m+1 \ (m\geq 1),
\end{cases}
\end{align*}
which is a proper discretisation of 
$\px ^{-1}(\frac{v^2}{2}) = \int _{-L}^x (\frac{v^2}{2})\, \mathrm dx$,
and satisfies
\begin{align*}
\delta _x \delta _x^{-1} \bra{\dfrac{1}{2}(v^n)^2} = \dfrac{1}{2}(v^n)^2
\end{align*}
for $n=1,\dots,N$.

Next we define the additional condition $v^0 =0$,
which corresponds to $u_x(-L,t)=0$, and calculate the temporal differentiation of $\mathcal H_{2,\rm d}$.
To simplify the calculation, we apply $\delta _x$ to \eqref{h2sds}
to obtain
\begin{align}
\delta _x \dot{u}^{n} = - \delta _x (u^nv^n) + \frac{1}{2}(v^n)^2. \label{h2sds2}
\end{align}
This relation holds for $n=1,\dots,N-1$.
We can now calculate $\dfrac{\mathrm d}{\mathrm dt} \mathcal H_{2,\rm d} (\bm u(t),\bm v(t))$, which directly leads to the following result.
In order to simplify the computations, we also assume $u^{-1} = u^1$, $v^{-1} = v^1$
and $u^{N+1}v^{N+1} = u^Nv^N$ so that \eqref{h2sds2} also holds even for $n=0,N$.
Then the temporal differentiation of $\mathcal H_{2,\rm d}$ is calculated as follows. 
\begin{align}
\dfrac{\mathrm d}{\mathrm dt} \mathcal H_{2,\rm d} (\bm u(t),\bm v(t))&= 
\dsum \Delta x \bra{\dfrac{(v^n)^2}{2} \dot{u}^n + u^n v^n \dot{v}^n} \\
&=
\dsum \Delta x \bra{\dfrac{(v^n)^2}{2} - \delta _x ( u^n v^n)}\dot{u}^n \\
&\quad+ \left[ \dfrac{1}{4} \bra{u^nv^n (\dot{u}^{n+1}+\dot{u}^{n-1}) + (u^{n+1}v^{n+1}+u^{n-1}v^{n-1}) \dot{u}^n} \right] _{0}^N \\
&=
\dsum \Delta x (\delta _x \dot{u}^n )\dot{u}^n 
+ \dfrac{1}{4} (u^{N+1}v^{N+1}+u^{N-1}v^{N-1}) \dot{u}^N \\
&= 
\left[ \dfrac{1}{4} (\dot{u}^{n+1} + \dot{u}^{n-1}) \dot{u^n} \right] _{0}^N + \dfrac{1}{4} (u^{N+1}v^{N+1}+u^{N-1}v^{N-1}) \dot{u}^N \\
&= \dfrac{1}{2} \dot{u}^{N-1}\dot{u}^N + \dfrac{1}{4} (u^{N+1}v^{N+1}+u^{N-1}v^{N-1}) \dot{u}^N \\
&= \dfrac{1}{2} \dot{u}^{N-1}\dot{u}^N + \dfrac{1}{2} u^{N-1}v^{N-1} \dot{u}^N. \label{hs:dpro}
\end{align}
Note that the discrete property \eqref{hs:dpro} is in good agreement with the continuous 
one \eqref{H2_td}, since $v^{N-1}$ (in the last term) is close to $0$ as will soon be observed in the numerical experiments.
 Applying the discrete gradient method for the temporal discretisation, we obtain the following fully-discrete scheme:
\begin{equation}\label{h2fds}
\delta _t^+ u^{n,i} = - u^{n,i+1/2} v^{n,i+1/2} + \delta _x^{-1} \bra{\dfrac{(v^{n,i+1})^2+(v^{n,i})^2}{4}} 
\end{equation}
for $n=1,\dots,N$,
still with the numerical boundary conditions
$u^{0,i} = 0$ and $v^{N,i} = 0$,
which correspond to $u(-L,t)=0$ and $u_x(L,t)=0$, respectively, for all $i$.
By construction, the scheme \eqref{h2fds} has the following property:
\begin{align*}
&\dfrac{1}{\Delta  t} ( \mathcal H_{2,\rm d} (\bm u^{i+1} , \bm v^{i+1}) - \mathcal H_{2,\rm d} (\bm u^{i} , \bm v^{i}) ) \\
&=
\dfrac{1}{2} (\delta _t^+ u^{N-1,i}) (\delta _t^+ u^{N,i})
+ \dfrac{1}{2} u^{N-1,i+1/2} v^{N-1,i+1/2} (\delta _t^+ u^{N,i}). 
\end{align*}

\subsubsection{Numerical simulations: Hamiltonian-preserving schemes for HS on $\hfline$}\label{ns_hsEP}
We now test the above presented  Hamiltonian-preserving integrators 
on the same problem as in Subsection~\ref{ns_hsMS}.

\begin{figure}
\centering
\includegraphics*[height=4.3cm,keepaspectratio]{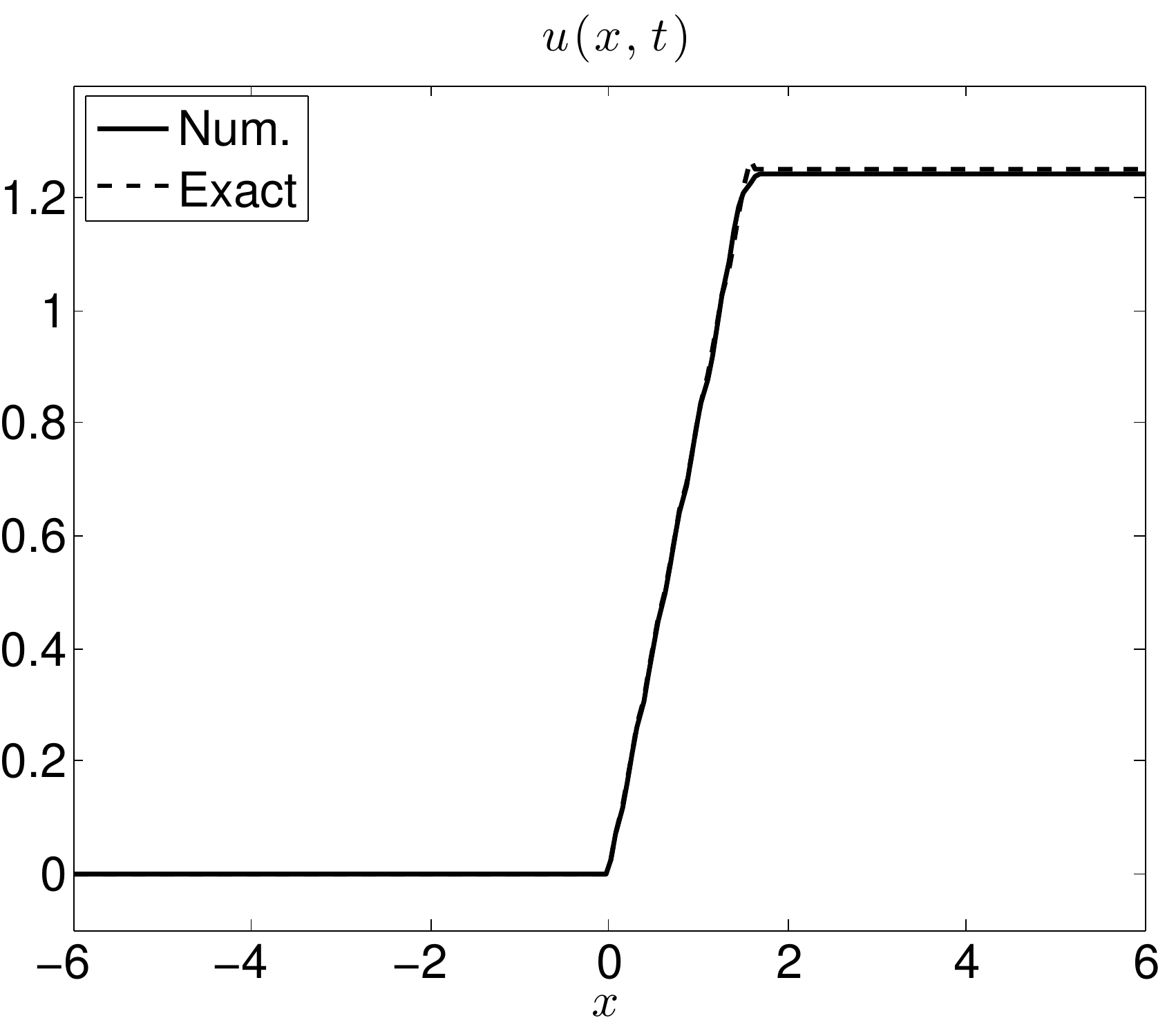}
\includegraphics*[height=4.3cm,keepaspectratio]{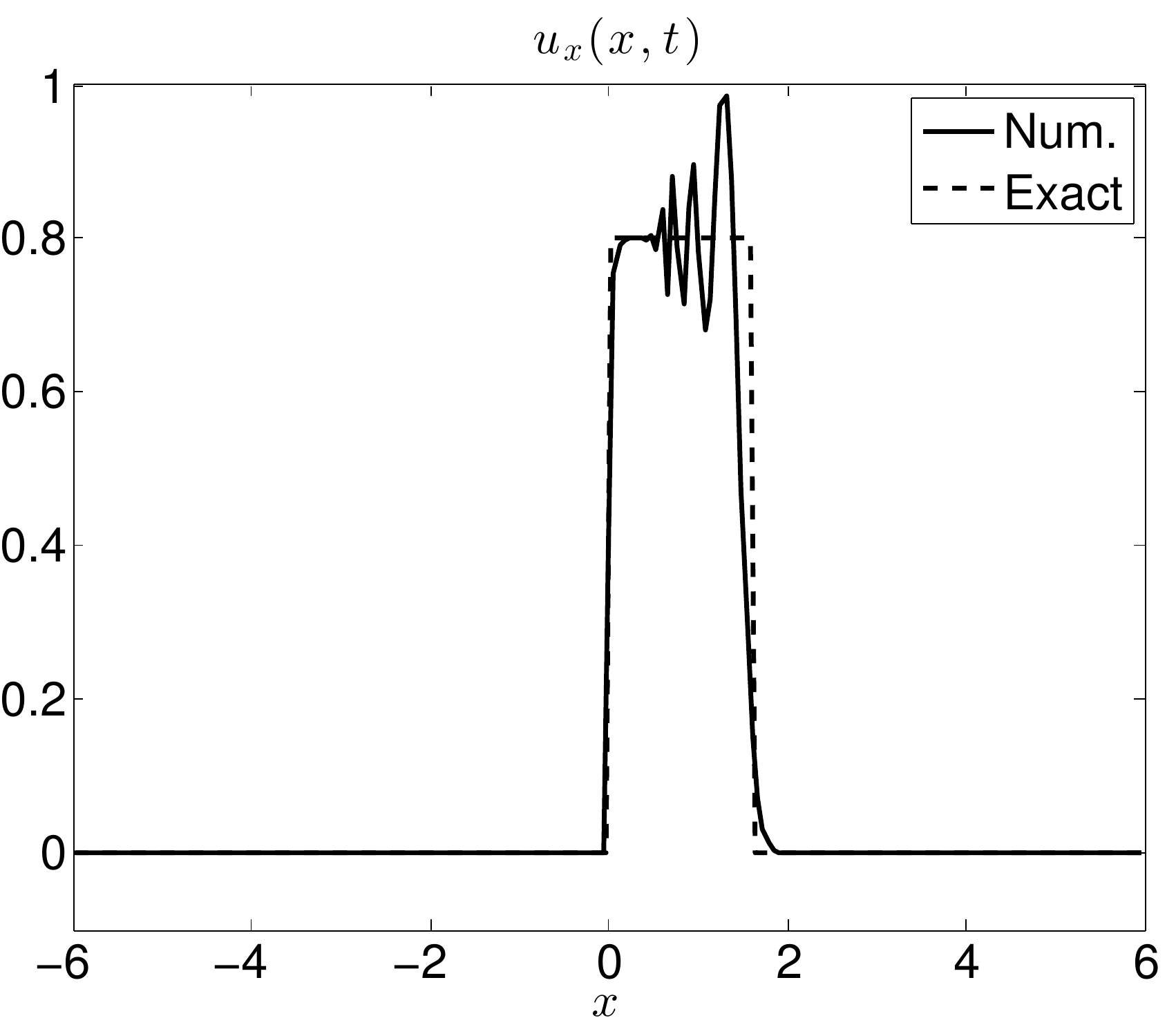}
\includegraphics*[height=4.3cm,keepaspectratio]{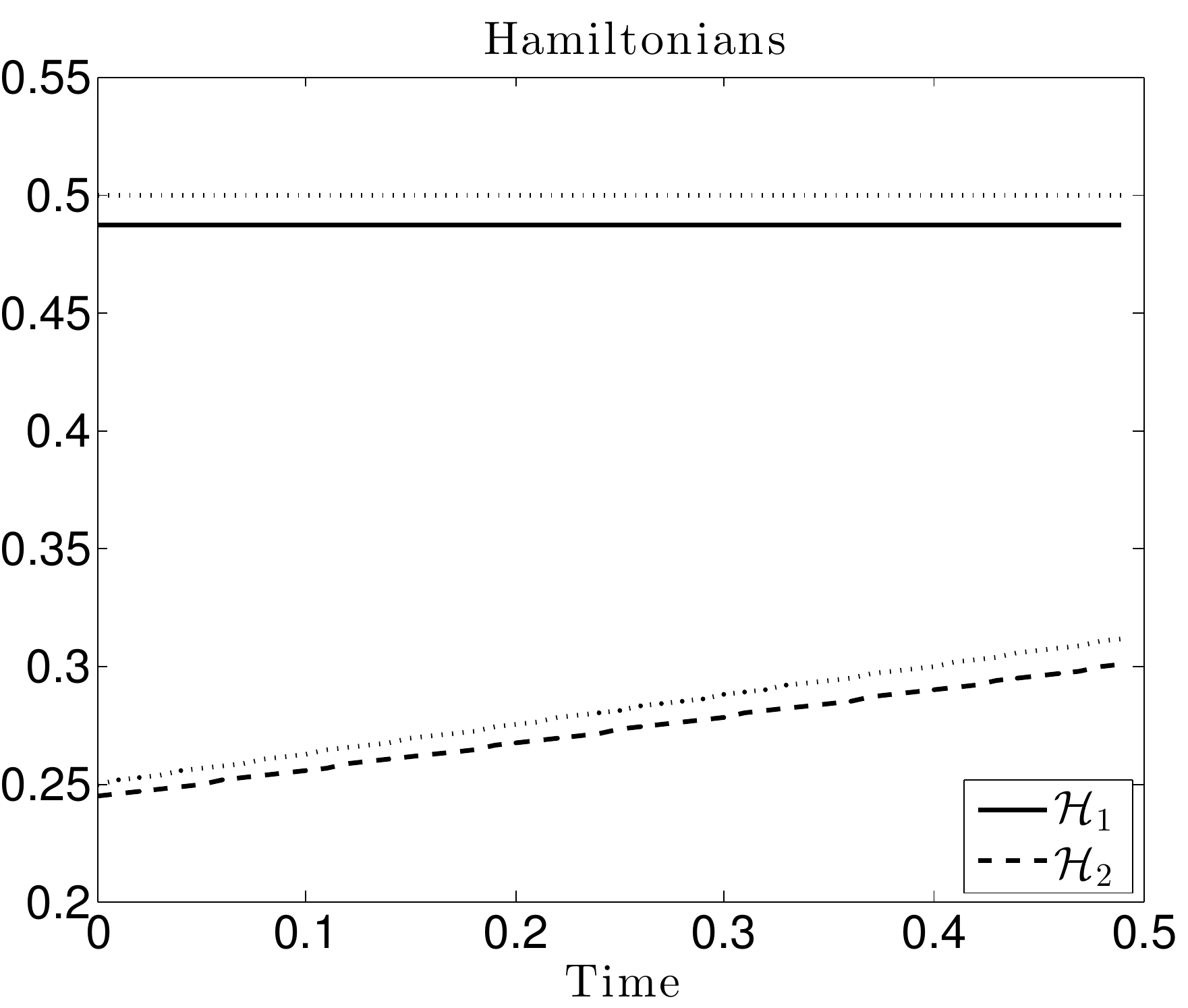}
\caption{$\mathcal H_1$-preserving scheme: exact and numerical profiles of $u(x,t)$ and $u_x(x,t)$ at time $T_{\text{end}} = 0.5$
and computed Hamiltonians ($\Delta x = 12/201$ and $\Delta t=0.01$).}
\label{fig:h1}
\end{figure}

\begin{figure}
\centering
\includegraphics*[height=4.3cm,keepaspectratio]{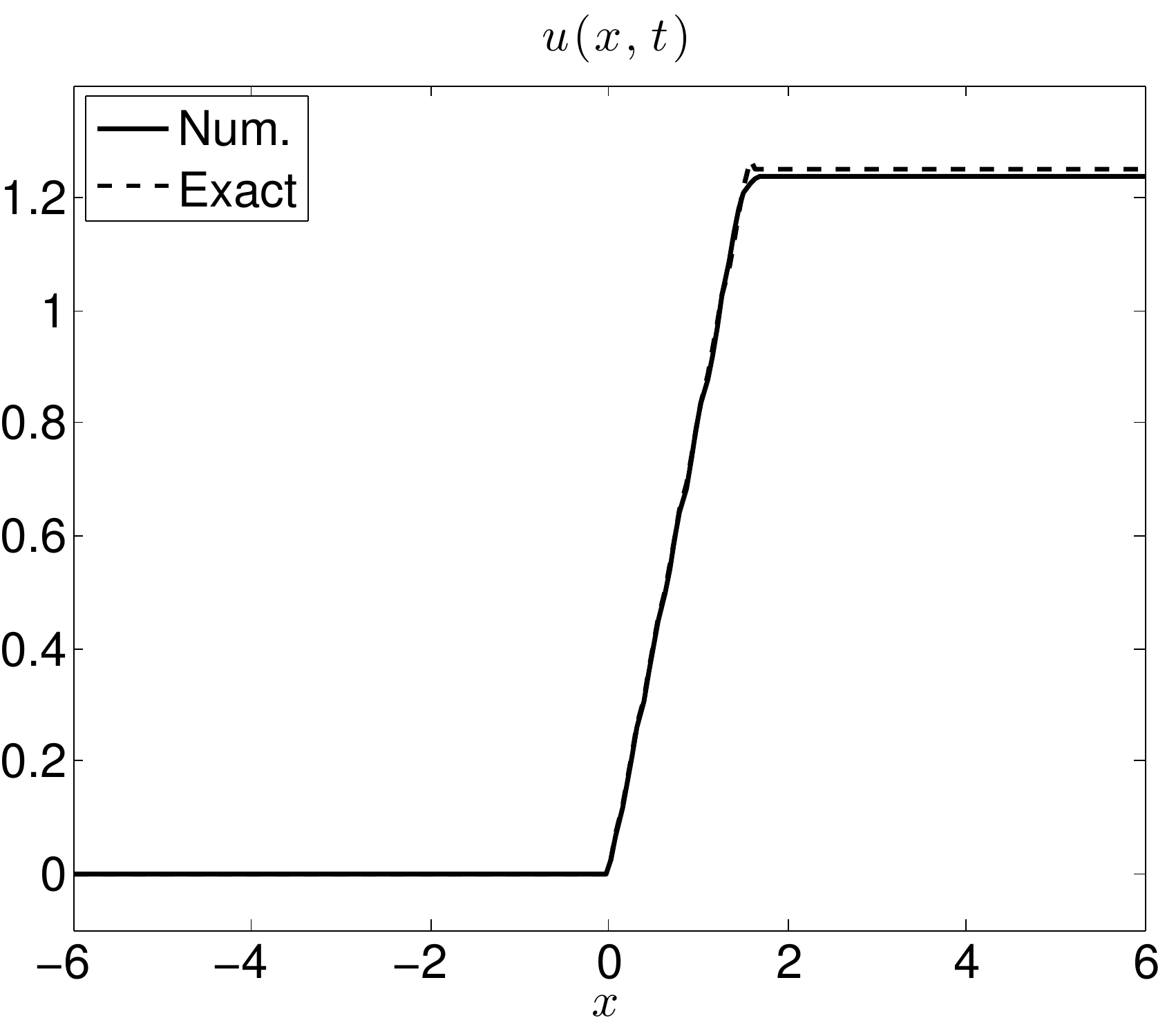}
\includegraphics*[height=4.3cm,keepaspectratio]{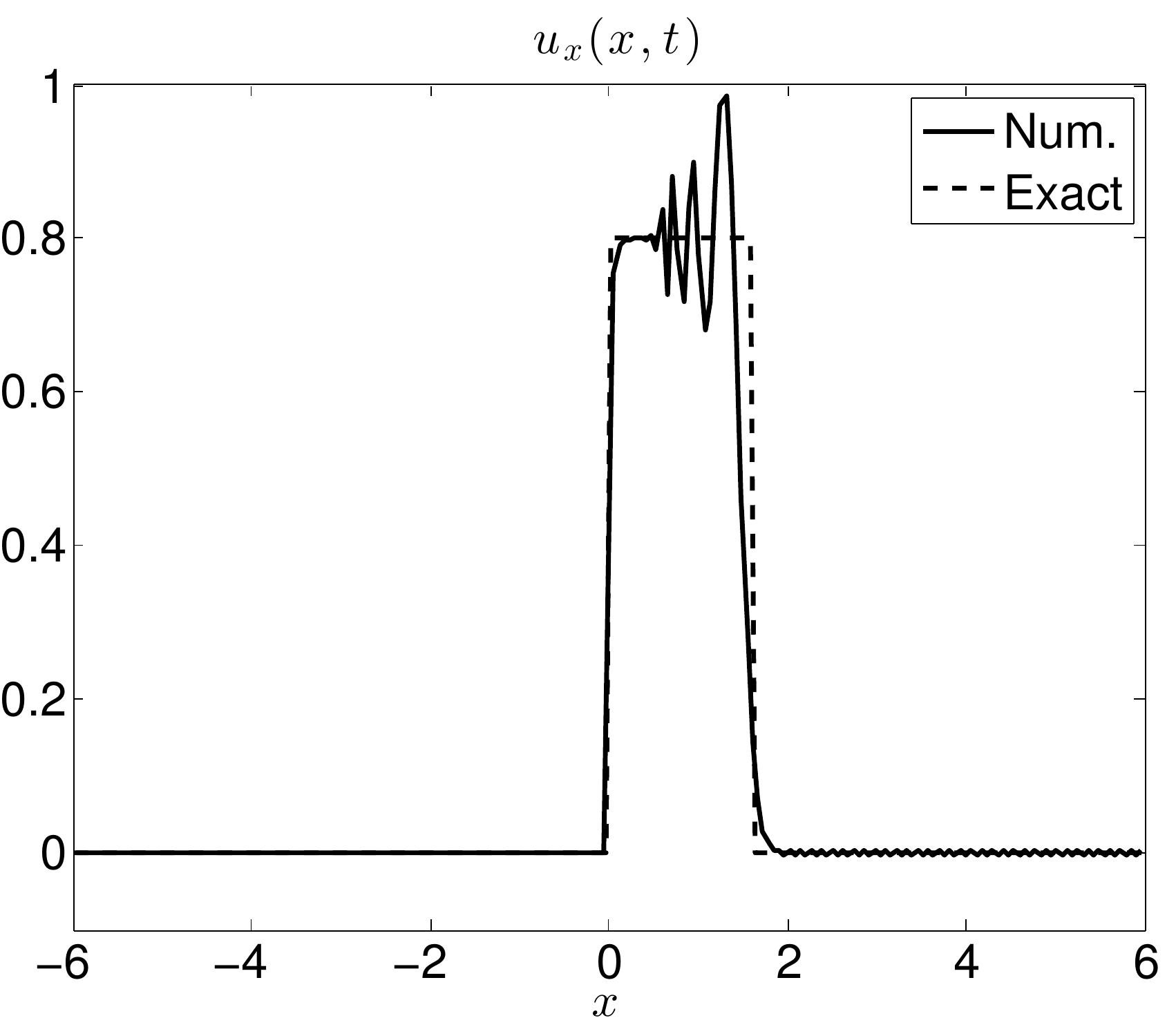}
\includegraphics*[height=4.3cm,keepaspectratio]{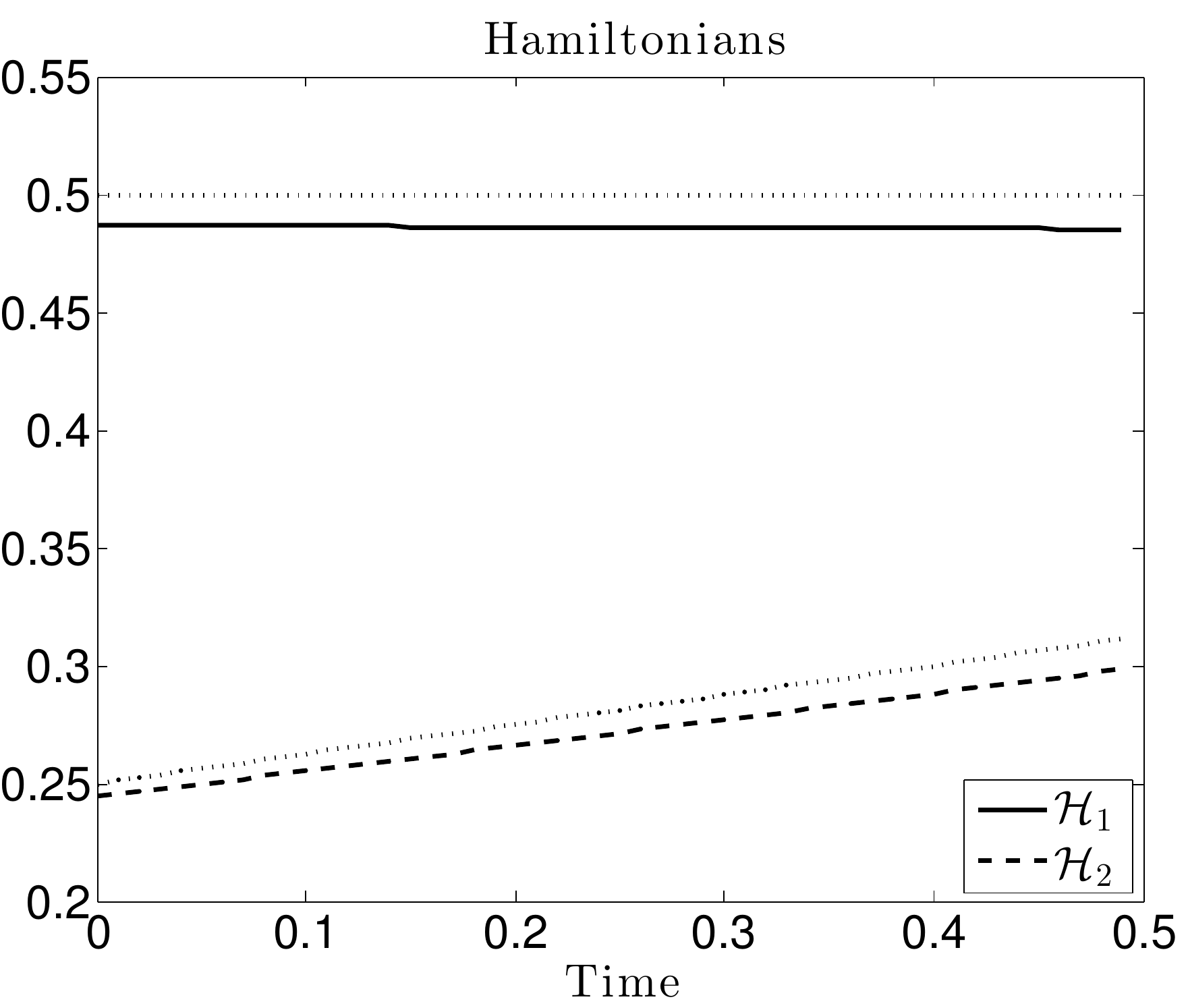}
\caption{$\mathcal H_2$-preserving scheme: exact and numerical profiles of $u(x,t)$ and $u_x(x,t)$ at time $T_{\text{end}} = 0.5$
and computed Hamiltonians ($\Delta x = 12/201$ and $\Delta t=0.01$).}
\label{fig:h2_2}
\end{figure}

Figures~\ref{fig:h1} and \ref{fig:h2_2} show the numerical results obtained by 
the $\mathcal H_1$-preserving scheme \eqref{h1fds} 
and $\mathcal H_2$-preserving scheme \eqref{h2fds} with step sizes 
$\Delta x=12/201$ (i.e., $N=201$) and $\Delta t=0.01$. 
Despite the oscillations observed in $u_x(x,t)$, one can observe 
a correct numerical approximation of $u(x,t)$.
Observe that the $\mathcal H_1$-preserving scheme offers a good behaviour in terms of the $\mathcal H_2$-quantity. 
For the $\mathcal H_2$-preserving scheme, though the slope of $\mathcal H_2$ does not perfectly 
coincide with that of the exact solution, the numerical result almost coincides with the one 
from the $\mathcal H_1$-preserving scheme. Furthermore, we also observe 
a good behaviour in terms of the $\mathcal H_1$-quantity. 
From a viewpoint of qualitative behaviour, no notable differences of 
the numerical solutions are observed 
for these two different energy-preserving schemes.

\subsection{Hamiltonian-preserving integrator for the modified Hunter--Saxton equation (periodic case)}
In this subsection, we propose an $\mathcal H_1$-preserving finite difference scheme for the mHS equation.
In order to derive the scheme,  
we use the first Hamiltonian structure~\cite{l08b}:
\begin{align*}
u_{xxt} = \left( (u_{xx}-\omega)\px + \px (u_{xx}-\omega ) \right) \px ^{-2} \fracdel{\mathcal H_1}{u},
\quad 
\mathcal{H}_1(u,u_x)=\frac12\int_0^L  u_{x}^2\,\mathrm{d}x,
\end{align*}
which can be expressed more explicitly as 
\begin{align*}
u_{xxt} = -\left( u_{xx}\px + \px u_{xx} \right) u + 2\omega u_x.
\end{align*}
In fact, one can prove the $\mathcal H_1$-preservation \eqref{hamilmhs} as follows:
\begin{align}
\label{mHScl1}
\dfrac{\mathrm d}{\mathrm dt} \mathcal H_1 
&= \int _0^L u_xu_{xt}\, \mathrm dx = -\int_0^L uu_{xxt} \,\mathrm dx
= \int _0^L u\cdot (u_{xx}\px + \px u_{xx})u \, \mathrm dx -2\omega \int _0^L uu_x \, \mathrm dx = 0. 
\end{align}
Here, no boundary terms appear due to the periodicity of solutions.
The last equality follows from the skew-symmetry of $(u_{xx}\px + \px u_{xx})$ and $\px$.

Based on this formulation, one derives the following $\mathcal H_1$-preserving scheme
\begin{equation}\label{ph1modHS2}
\tilde{\delta}_x^2 \delta_t^+ u^{n,i} + (\tilde{\delta}_x^2 u^{n,i+1/2}) (\delta_x u^{n,i+1/2}) + \delta_x (u^{n,i+1/2} (\tilde{\delta}_x^2 u^{n,i+1/2}))
- 2\omega \delta _x u^{n,i+1/2} = 0,
\end{equation}
for which the solutions satisfy
\begin{align*}
\mathcal H_{1,\rm d} (\bm{u}^i) = \mathcal H_{1,\rm d} (\bm{u}^0) \quad
{\rm for \ all} \ i,
\end{align*}
where \begin{align*}
\mathcal H_{1,\rm d} (\bm{u}^i) := {\sum_{n=0}^{N-1}}   \Delta x \dfrac{(\delta_ x^+ u^n)^2}{2}.
\end{align*}
This conservation property is checked by calculations corresponding to \eqref{mHScl1} in the discrete setting.

As this was done for the multi-symplectic integrator,
in order to select the right travelling waves, we have to consider 
the pseudo-inverse operator of $\tilde{\delta} _x^2$, denoted by $(\tilde{\delta} _x^2)^\dagger$. 
We then obtain the following scheme
\begin{equation}\label{ph1modHS}
\delta_t^+ u^{n,i} + (\tilde{\delta} _x^2) ^\dagger \left((\tilde{\delta}_x^2 u^{n,i+1/2}) (\delta_x u^{n,i+1/2}) + \delta_x (u^{n,i+1/2} (\tilde{\delta}_x^2 u^{n,i+1/2}))
- 2\omega \delta _x u^{n,i+1/2} \right) = 0.
\end{equation}
The above conservation property remains
since the solution of \eqref{ph1modHS} always satisfies \eqref{ph1modHS2}.

\begin{figure}%
\begin{center}
\includegraphics*[scale=0.3]{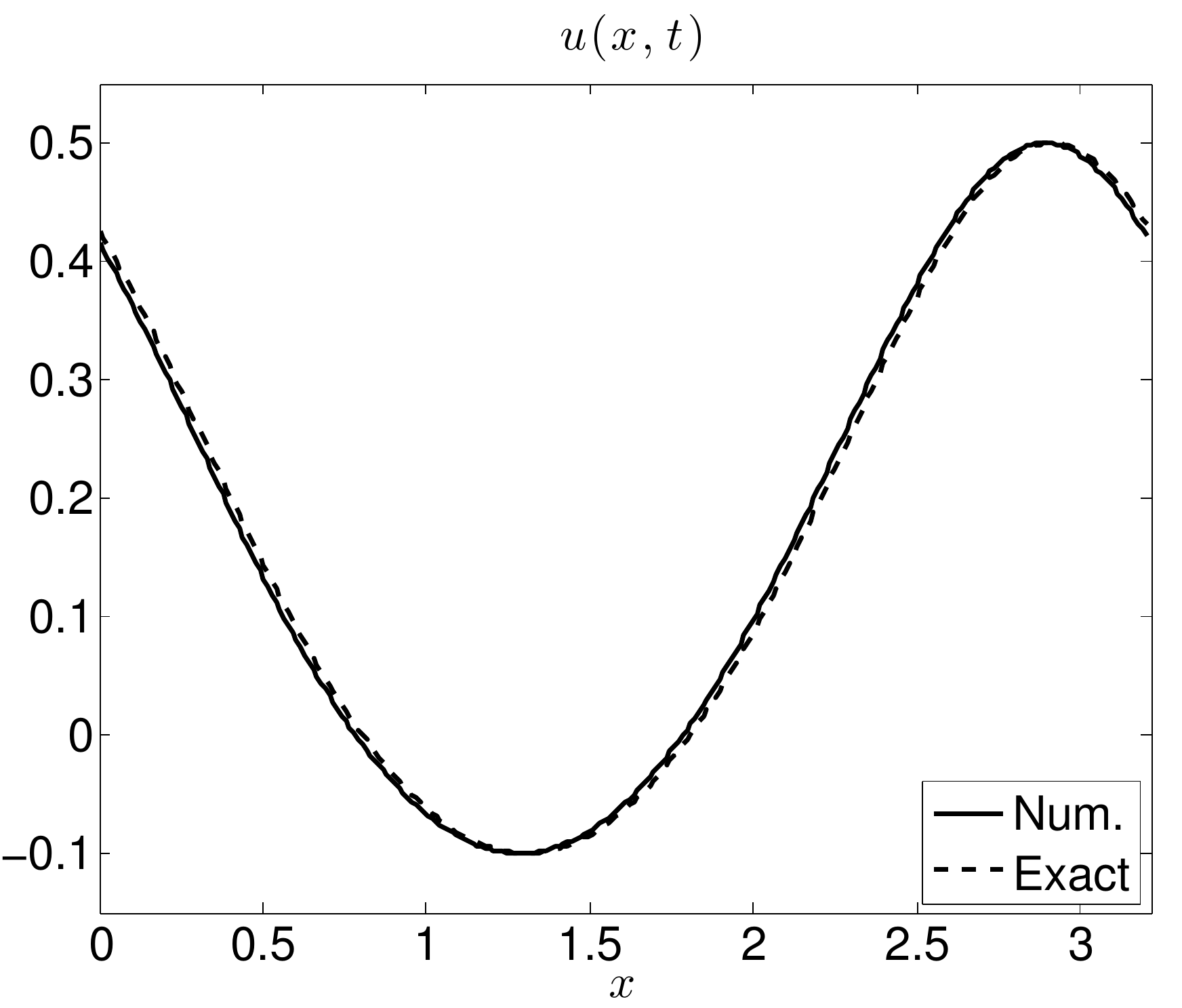}
\includegraphics*[scale=0.3]{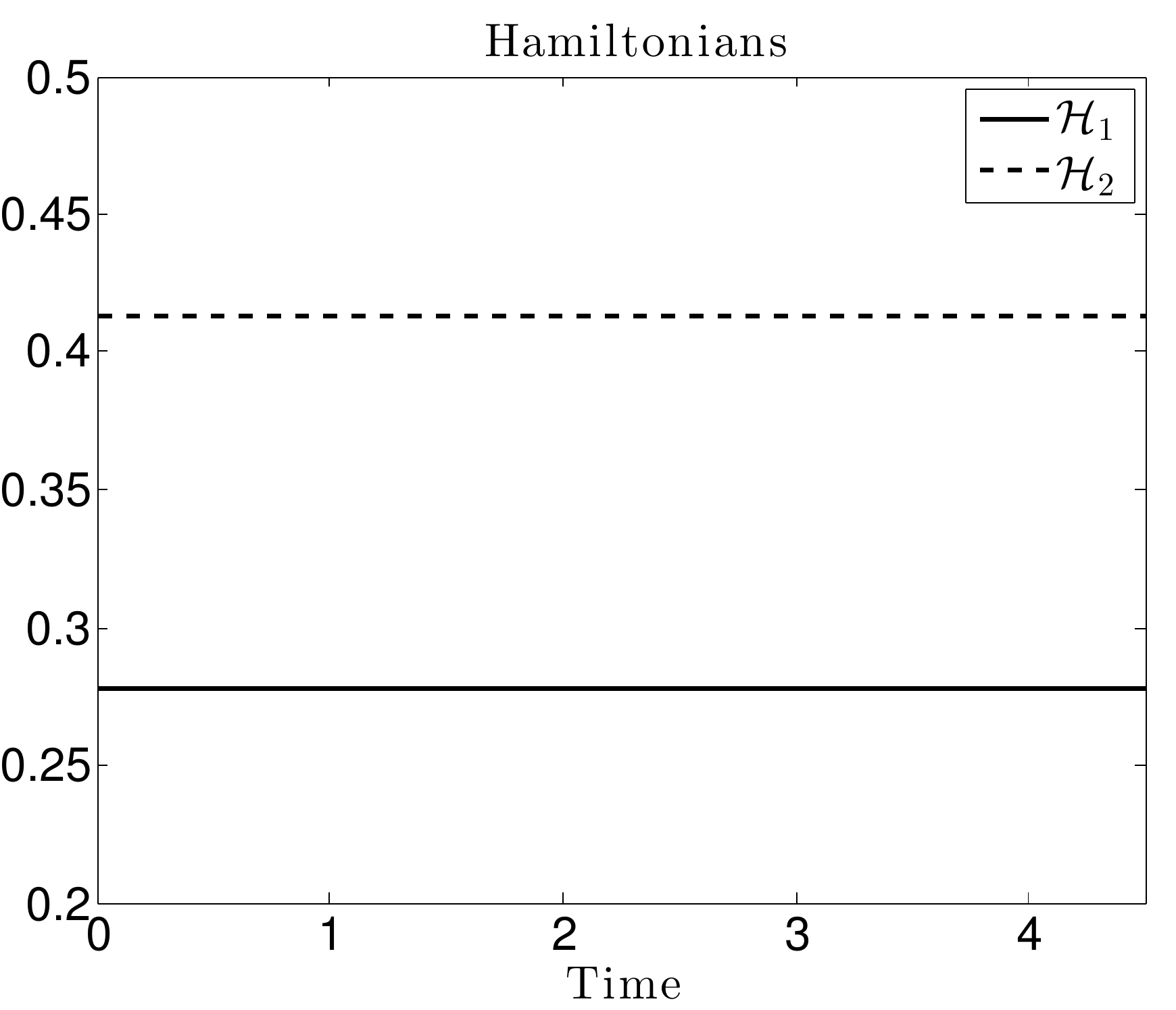}
\caption{Exact and numerical profiles of $u(x,t)$  
at time $T_{\text{end}}=3.5$ and computed Hamiltonians by the $\mathcal H_1$-preserving scheme \eqref{ph1modHS}. 
}
\label{fig:h1mod}
\end{center}
\end{figure}

We now apply the scheme \eqref{ph1modHS} to the same problem as in Subsection~\ref{mssimmHS}. 
Figure~\ref{fig:h1mod} displays the exact and 
numerical profiles of $u(x,t)$ at time $T_{\text{end}}=3.5$ 
and also the computed values of the Hamiltonians 
\begin{align*}
\mathcal{H}_1(u,u_x)=\frac12\int_0^L  u_{x}^2\,\mathrm{d}x,\quad
\mathcal{H}_2(u,u_x)=\frac12\int_0^L (uu_x^2+2\omega u^2)\,\mathrm{d}x
\end{align*} 
with (relative large) step sizes 
$\Delta t=0.02$ and $\Delta x=L_{\text{per}}/256$. 

One may notice that the numerical solutions given by this integrator agree very well 
with the exact ones as well as with the ones given by the multi-symplectic schemes from 
the previous section. In addition to the exact preservation of $\mathcal{H}_1$, the quantity $\mathcal{H}_2$ is also 
preserved with a very good accuracy.

\subsection{Hamiltonian-preserving integrator for the two component Hunter--Saxton equation (periodic case)}

In this subsection, we propose an $\mathcal H_1$-preserving finite difference scheme for the 2HS equation.
In order to derive the scheme,  
we use the first Hamiltonian structure~\cite{l09}:
\begin{align*}
\begin{bmatrix}
u_{xxt}\\ \rho_t
\end{bmatrix}=
\begin{bmatrix}
(u_{xx}\px + \px u_{xx}) & \rho \px \\ \px \rho & 0
\end{bmatrix}
\begin{bmatrix}
\px ^{-2} \delta \mathcal H_1 / \delta u \\ \delta \mathcal H_1 / \delta \rho
\end{bmatrix}, \quad
\mathcal{H}_1(u,u_x,\rho)=\frac12\int(u_x^2+\kappa\rho^2)\,\mathrm{d}x,
\end{align*}
which can be expressed more explicitly as 
\begin{align*}
& u_{xxt} = -\left( u_{xx}\px + \px u_{xx} \right) u + \kappa \rho \rho _x,\\
& \rho _t = - (u\rho )_x .
\end{align*}
In fact, one can prove the $\mathcal H_1$-preservation \eqref{hamil2HS} as follows:
\begin{align*}
\dfrac{\mathrm d}{\mathrm dt} \mathcal H_1 
&= \int _0^L u_xu_{xt}\, \mathrm dx + \kappa \int _0^L \rho\rho _t\, \mathrm dx
= - \int _0^L uu_{xxt}\, \mathrm dx - \kappa \int _0^L \rho (u\rho )_x\, \mathrm dx \nonumber \\ 
&= \int _0^L u\cdot (u_{xx}\px + \px u_{xx}) u\, \mathrm dx - \kappa \int _0^L u \rho \rho _x\, \mathrm dx
- \kappa \int _0^L \rho (u\rho )_x\, \mathrm dx \nonumber \\
&= - \kappa \int _0^L u \rho \rho _x \,\mathrm dx + \kappa \int _0^L \rho_x u\rho\, \mathrm dx = 0. 
\end{align*}
As it was done previously, one derives the following $\mathcal H_1$-preserving scheme
\begin{eqnarray}
\begin{array}{l}
\delta_t^+ u^{n,i} + (\tilde{\delta} _x^2) ^ \dagger \left((\tilde{\delta}_x^2 u^{n,i+1/2}) (\delta_x u^{n,i+1/2}) + \delta_x (u^{n,i+1/2} (\tilde{\delta}_x^2 u^{n,i+1/2}))
- \kappa \rho^{n,i+1/2}  \delta _x \rho^{n,i+1/2} \right) = 0, \\
\delta _t ^+ \rho ^{n,i} + \delta _x (u^{n,i+1/2}\rho^{n,i+1/2}) = 0,
\end{array}
\label{ph12HS}
\end{eqnarray}
with the pseudo-inverse operator $(\tilde{\delta} _x ^2) ^\dagger$.
This scheme possesses a conservation property
\begin{align*}
\mathcal H_{1,\rm d} (\bm{u}^i) = \mathcal H_{1,\rm d} (\bm{u}^0) \quad
{\rm for \ all} \ i,
\end{align*}
where \begin{align*}
\mathcal H_{1,\rm d} (\bm{u}^i) := {\sum_{n=0}^{N-1}}   \Delta x \dfrac{(\delta_ x^+ u^n)^2 + \kappa (\rho ^n)^2}{2} .
\end{align*}

We now test the $\mathcal H_1$-preserving scheme \eqref{ph12HS} on the 
travelling wave solutions of \eqref{eq:HSsyst} with $\kappa=1$ 
used in Subsection~\ref{mssim2HS}. 

\begin{figure}%
\begin{center}
\includegraphics*[height=4.4cm,keepaspectratio]{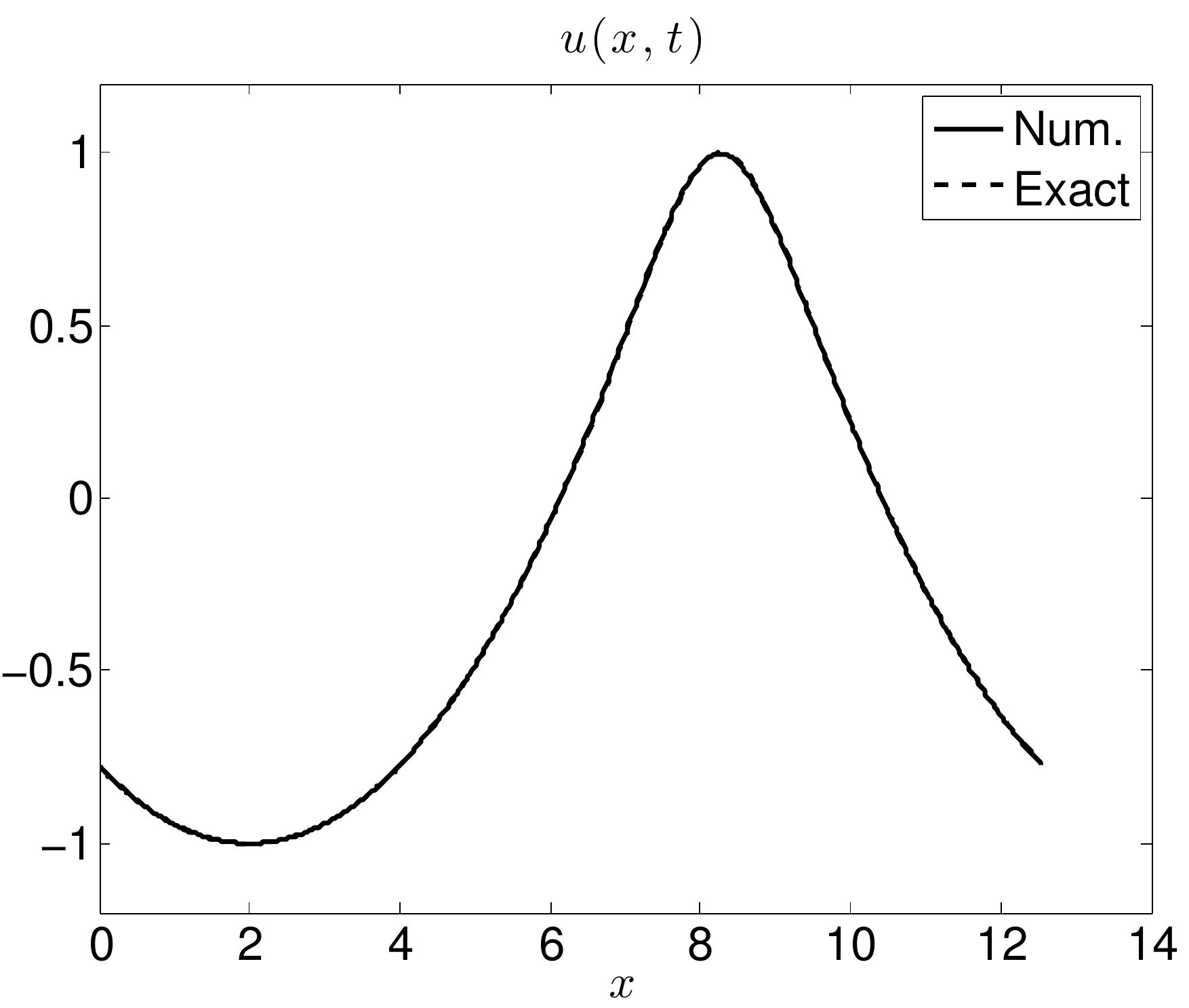}
\includegraphics*[height=4.4cm,keepaspectratio]{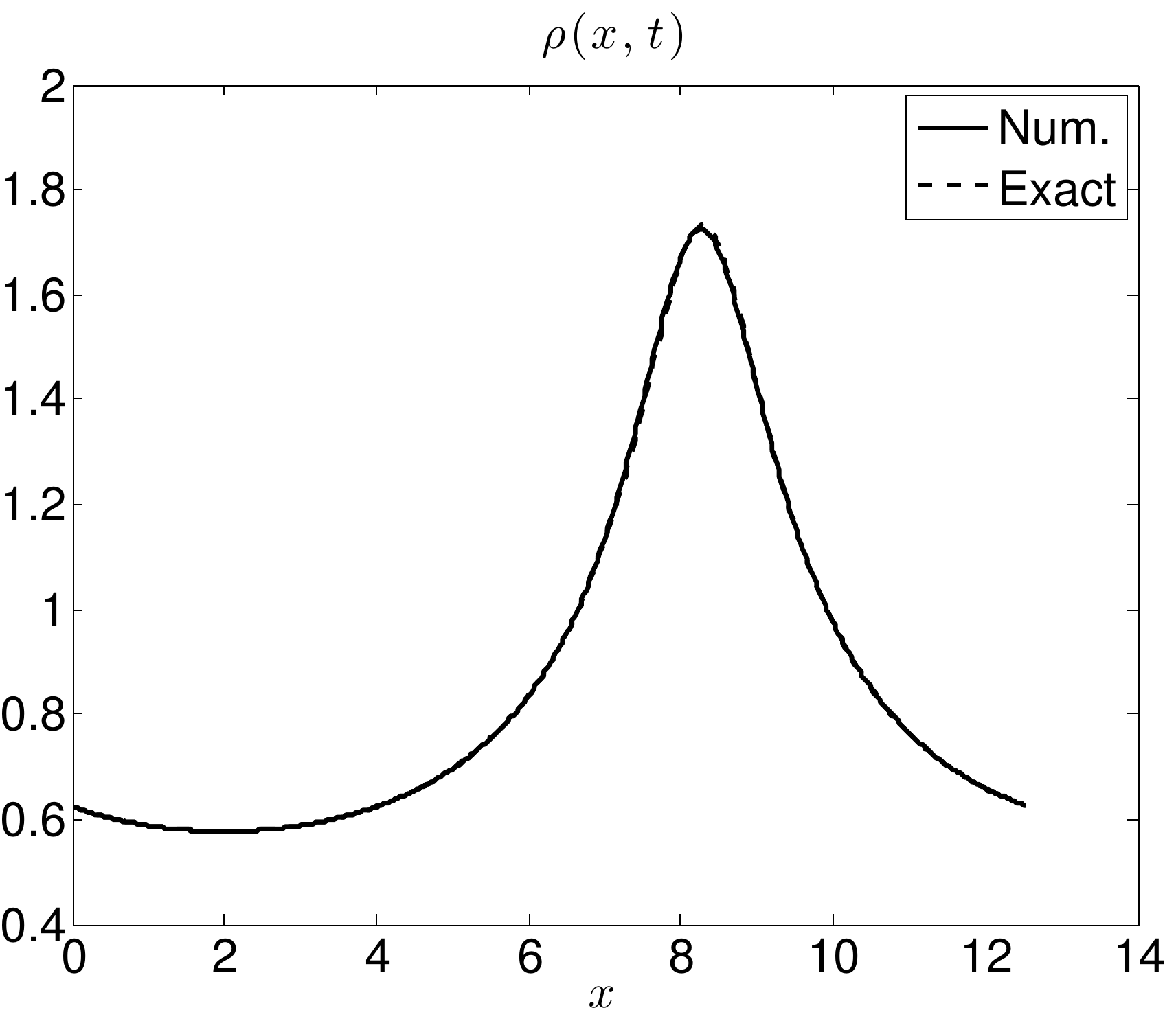}
\includegraphics*[height=4.4cm,keepaspectratio]{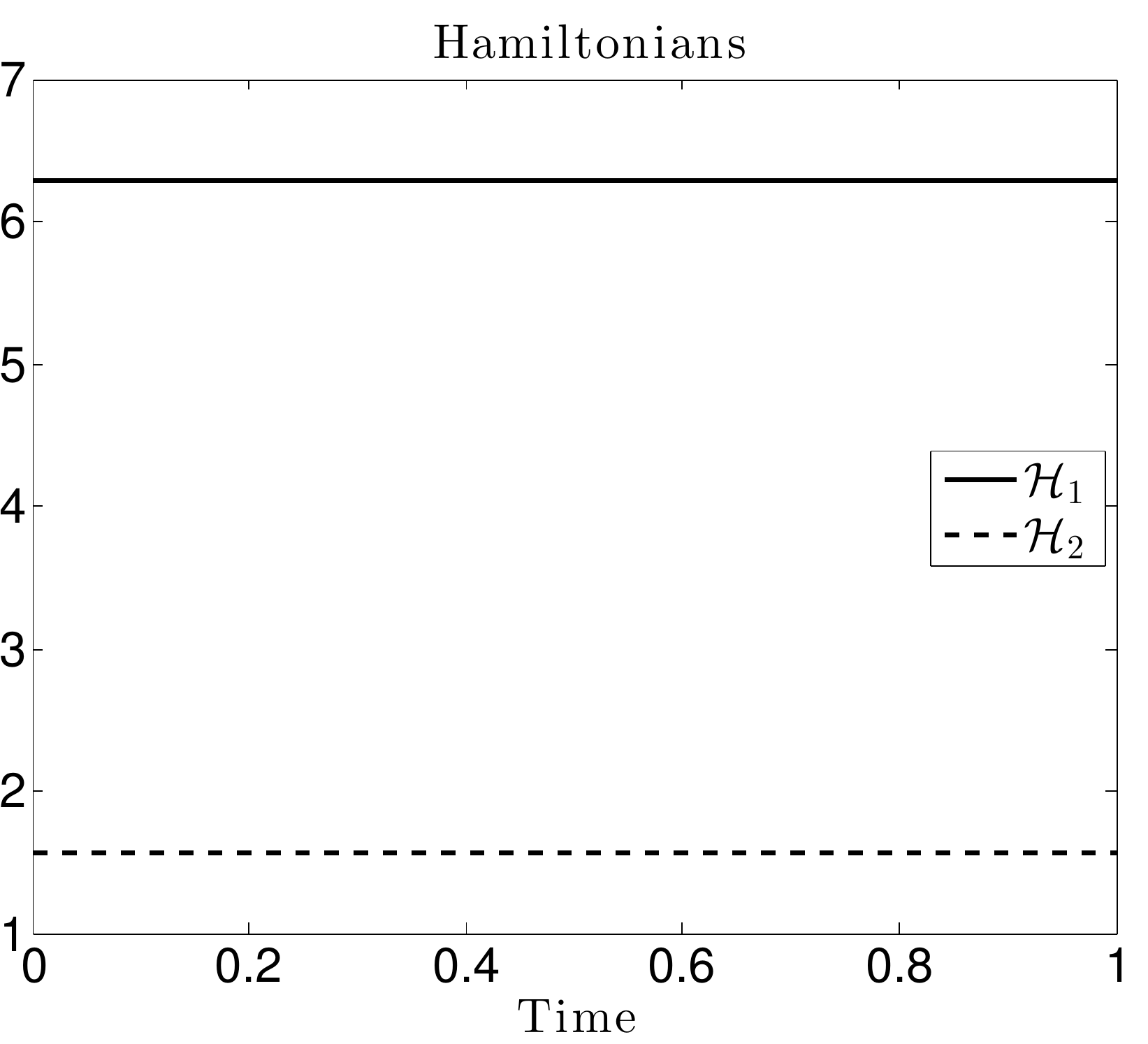}
\caption{Exact and numerical profiles of $u(x,t)$ and $\rho(x,t)$   
at time $T_{\text{end}}=1$ and computed Hamiltonians (right plot)
by the $\mathcal H_1$-preserving scheme \eqref{ph12HS}. 
}
\label{fig:cons2}
\end{center}
\end{figure}

Figure~\ref{fig:cons2} displays the exact and 
numerical profiles of $u(x,t)$ and $\rho(x,t)$ at time $T_{\text{end}}=1$ 
and also the computed values of the Hamiltonians \eqref{hamil2HS} 
using the $\mathcal H_1$-preserving scheme \eqref{ph12HS} 
with (large) step sizes 
$\Delta t=0.1$ and $\Delta x=L_{\text{per}}/512$.

The numerical profiles of the solutions agree very well with the exact profiles as well as 
with the profiles of the multi-symplectic scheme displayed in Figure~\ref{fig:trav2}. 
Further, excellent conservation properties of the numerical solutions are observed. 

\section{Concluding remarks}\label{sec-conc}
This study contributes to enhance multi-symplectic discretisations 
for partial differential equations arising from important applications 
in physics. 
Indeed, we have presented the first multi-symplectic 
formulations of the Hunter--Saxton equation, 
of the modified Hunter--Saxton equation and 
of the two component Hunter--Saxton system. 
Furthermore, using these results, we have derived 
novel explicit multi-symplectic integrators for these nonlinear 
partial differential equations. 

As exact solutions to these PDEs are rarely known, further numerical 
methods are useful for comparison. We therefore investigate Hamiltonian-preserving 
numerical discretisations for these problems. No significant differences between this  
type of geometric numerical schemes and the multi-symplectic integrators are observed 
except, perhaps, the fact that the last one seem 
a little bit better from a practical viewpoint, being explicit and a little 
easier to implement. 

A major difficulty in the derivation of numerical schemes for the HS-like equations is 
a proper treatment of the boundary conditions. We clarify this issue for the exact as well as for 
the numerical solutions to these PDEs. Therefore, all the numerical methods proposed in this publication enjoy a correct treatment of the boundary conditions for the HS-like equations.

So far, numerical tests have been conducted only 
with the Euler box scheme. Besides numerical experiments have 
been conducted for particular travelling wave solutions 
of the modified Hunter--Saxton equation or the two component 
Hunter--Saxton system. It thus remains to try out and 
analyse more elaborate structure-preserving 
numerical methods and other type of solutions to these problems.

\bibliographystyle{plain}
\bibliography{bibhs}

\end{document}